%%%%%%%%%%%%%%%%%%   Geometry and Topology: 2000-12.tex  %%%%%%%%%%%%%
%%%%        
%%%%           Taut ideal triangulations of 3-manifolds
%%%%             
%%%%                           Marc Lackenby
%%%%  
%%%%             Published in Volume 4(2000) pages 369-395
%%%%
%%%%                   Publication date 4 November 2000
%%%%
%%%%                      This is a plain TeX file
%%%%
%%%%
%%%%%%%%%%%%%%%%%%                                   %%%%%%%%%%%%%%%%%%%

%%%%%%%%%%%%%%%%%%%%%%%%%%%%%%%%%%%%%%%%%%%%%%%%%%%%%%%%%%%%%%%
%%%%%%%%%%%             gtmacros.tex            %%%%%%%%%%%%%%%
%%%%%%%%%%%             version 1.6             %%%%%%%%%%%%%%% 
%
%                       Colin Rourke   
%
%
%    These macros are recommended for use by authors submitting articles   
%    to Geometry and Topology or to Algebraic and Geometric Topology.  
%    They are intended to be used with plain TeX. Each macro is described 
%    briefly to make it clear how to use it (or to modify it to achieve
%    different results).  If you modify this file then please change its
%    name.  If you modify this file and use the modified file to 
%    format an article for submission to Geometry and Topology or
%    Algebraic and Geometric Topology, then please paste the modified
%    file into your main TeX file.  Do not submit it as a separate file.
%      
%    Instructions on using these macros are also given in  gtmacins.tex  
%    or  gtmacins.ps  or .pdf  available on the gt www pages or by 
%    anonymous ftp from the gt/info/macros directory.
%
%
\magnification=\magstephalf      % Sets default point size to 11pt.
%
%  Basic layout parameters :
%
\vsize=7.5truein                 % Sets text height to 7.5 inches.
\hsize=5.2truein                 % Sets text width to 5.2 inches.
\newskip\stdskip                 % standard vertical space
\stdskip=6pt plus3pt minus3pt    % (slightly more stretchy
\medskipamount=\stdskip          % than the usual \medskip)
\parindent=0pt                   % Paragraphs are non-indented with
\parskip=\stdskip                % a little space between paragraphs. 
\abovedisplayskip=\stdskip       %  Reduces the space
\belowdisplayskip=\stdskip       %  around displays.
\mathsurround=0.75pt             % Gives a little extra space around maths.
\overfullrule=0pt                %  Prevents black boxes
%
%   The following macro is for principal paragraph breaks ie
%   a paragraph break with a slightly larger space :
%
\def\ppar{\par\goodbreak\vskip 8pt plus 4pt minus 4pt}     
%
%  The standard horizontal space for theorems, labels etc :
%
\def\stdspace{\hskip 0.75em plus 0.15em\ignorespaces}
\let\qua\stdspace % useful abbreviation (3/4 of a quad)
%
%%%%%%%%%%%%%%            FONT MACROS            %%%%%%%%%%%%
%
%           The following font macros define the AMS symbol 
%           and Euler-Fraktal fonts for use in text and
%           mathematics with appropriate size changes.
%           They also define two new control sequences  
%           \small  and  \large  (similar to those built
%           into LaTeX) which change the size of all fonts 
%           both in text and maths.  \small  is 10% smaller 
%           than normal and  \large  30% bigger.  The strange
%           size of the \large text fonts (10pt scaled 1315)
%           is because these macros are intended to be used
%           at \magstephalf.  The result is 10pt scaled 1440
%           (\magstep2) which is a standard font size.  If
%           you are borrowing these macros to use them at
%           another basic  \magnification, then you will
%           probably need to change 1315 to 1200 in the eleven
%           places marked ** below.  \large  will then be
%           20% bigger than normal.  Note that at \magstephalf
%           all the fonts come out roughly one point larger
%           than their size as defined in these macros.
%
%           The size-changing macros are based on Knuth's
%           \ninepoint and \eightpoint macros.
%
%
%    The macros are laid out in a way which makes it clear how to
%    add futher fonts (or delete unavailable ones) and how to add
%    further size changes.
%
%    First comes a definition of  \hexnumber  which is needed to
%    refer to font families whose family number is not known :
%
\def\hexnumber#1{\ifcase#1 0\or 1\or 2\or 3\or 4\or 5\or 6\or 7\or 8\or
 9\or A\or B\or C\or D\or E\or F\fi}
%
%     Next we define the AMS symbol-a fonts at 13,10,9,7,6,5 points
%
\font\thirtnmsa=msam10 scaled 1315    %%% **  see note above 
\font\tenmsa=msam10          \font\ninemsa=msam9
\font\sevenmsa=msam7         \font\sixmsa=msam6
\font\fivemsa=msam5
%%%%%%  (add further sizes here if you need them)
%
%    and the standard size family for these fonts
%
\newfam\msafam                  \textfont\msafam=\tenmsa
\scriptfont\msafam=\sevenmsa    \scriptscriptfont\msafam=\fivemsa
\edef\hexa{\hexnumber\msafam}        %  The msa family is  \fam\hexa
\def\msa{\fam\msafam\tenmsa}         %  \msa  switches to this family
%
%    Repeat these steps for the AMS symbol-b fonts
%
\font\thirtnmsb=msbm10 scaled 1315   %%%  ** see note above
\font\tenmsb=msbm10      \font\ninemsb=msbm9
\font\sevenmsb=msbm7     \font\sixmsb=msbm6
\font\fivemsb=msbm5
%%%%%%  (add further sizes here if you need them)
%
\newfam\msbfam                   \textfont\msbfam=\tenmsb       
\scriptfont\msbfam=\sevenmsb     \scriptscriptfont\msbfam=\fivemsb
\edef\hexb{\hexnumber\msbfam}    %  The msb family is \fam\hexb  
\def\msb{\fam\msbfam\tenmsb}     %  \msb switches to this family
%
%        Repeat for the Euler-Fraktal fonts 
%
\font\thirtneufm=eufm10 scaled 1315   %%% **  see note above 
\font\teneufm=eufm10                 \font\nineeufm=eufm9
\font\seveneufm=eufm7                \font\sixeufm=eufm6
\font\fiveeufm=eufm5
%%%%%%  (add further sizes here if you need them)
%
\newfam\eufmfam                    \textfont\eufmfam=\teneufm
\scriptfont\eufmfam=\seveneufm     \scriptscriptfont\eufmfam=\fiveeufm
\edef\hexf{\hexnumber\eufmfam}      % The Euler-Fraktal family is
\def\frak{\fam\eufmfam\teneufm}     % \fam\hexf and \frak switches to this
%
%%%  Add further fonts families here (using the same format) if you need
%    them.  The def of hexnumber is optional (it is only used for
%    \mathchardef 's).
%
%      Now we need to define the standard fonts (which are
%      already defined at 10,7 and 5 point) at 13,9 and 6 point:
%
%      Roman fonts:
\font\thirtnrm=cmr10 scaled 1315    %%%  ** see note above
\font\ninerm=cmr9                   \font\sixrm=cmr6   
%%%%%%  (add further sizes here if you need them)
%
%      Math italic fonts
\font\thirtni=cmmi10 scaled 1315    %%%  ** see note above 
\font\ninei=cmmi9                   \font\sixi=cmmi6  
%%%%%%  (add further sizes here if you need them)
%
%     Symbol fonts
\font\thirtnsy=cmsy10 scaled 1315   %%%  ** see note above
\font\ninesy=cmsy9                  \font\sixsy=cmsy6  
%%%%%%  (add further sizes here if you need them)
%
%     Bold face
\font\thirtnbf=cmbx10 scaled 1315   %%%  ** see note above 
\font\ninebf=cmbx9                  \font\sixbf=cmbx6  
%%%%%%  (add further sizes here if you need them)
%
%     The maths extension font (only defined at text size)
%
\font\thirtnex=cmex10 scaled 1315   %%%  ** see note above
\font\nineex=cmex9                  
%%%%%%  (add further sizes here if you need them)
%
%     Finally three fonts (text italic, slanted and typewriter type)
%     which are also only defined at text size
%
\font\thirtnit=cmti10 scaled 1315  %%%  ** see note above 
\font\nineit=cmti9                  
%%%%%%  (add further sizes here if you need them)
%
\font\thirtnsl=cmsl10 scaled 1315  %%%  ** see note above 
\font\ninesl=cmsl9                  
%%%%%%  (add further sizes here if you need them)
%
\font\thirtntt=cmtt10 scaled 1315  %%%  ** see note above 
\font\ninett=cmtt9                  
%%%%%%  (add further sizes here if you need them)
%
%
%     Now come the two main macros.  What  \small  does is to
%     change all the families of fonts from normal size which is
%     10,7,5  (ie 10pt text, 7pt subscript, 5pt subsubscript)
%     to 9,6,5.  \large  similarly changes to  13,9,7.  To make
%     other size changing macros, choose your three sizes, add
%     font size definitions if necessary and make the obvious changes
%     to one of these macros.  Change  \normalbaselineskip  and
%     \strutbox  dimensions to appropriate sizes as well.  To
%     add further fonts, insert them in each macro, using the
%     AMS fonts as a model.
%      
%
\def\small{%
%
%   redefine the sizes of the roman fonts :
%
\textfont0=\ninerm \scriptfont0=\sixrm \scriptscriptfont0=\fiverm
\def\rm{\fam0\ninerm}%       % ( \rm  sets \ninerm  in text mode
%                            %  and \fam0 in math mode)
%
%   and the math italic fonts :
%
\textfont1=\ninei \scriptfont1=\sixi \scriptscriptfont1=\fivei
%
%   and the symbol fonts :
%
\textfont2=\ninesy \scriptfont2=\sixsy \scriptscriptfont2=\fivesy
%
%   There is only one math extension font :
%
\textfont3=\nineex \scriptfont3=\nineex \scriptscriptfont3=\nineex
%
%   Next the bold font (named rather than numbered) :
%
\textfont\bffam=\ninebf \scriptfont\bffam=\sixbf
\scriptscriptfont\bffam=\fivebf \def\bf{\fam\bffam\ninebf}%
%
%   and the three text-only fonts : 
%
\textfont\itfam=\nineit \def\it{\fam\itfam\nineit}%
\textfont\slfam=\ninesl \def\sl{\fam\slfam\ninesl}%
\textfont\ttfam=\ninett \def\tt{\fam\ttfam\ninett}%
%
%   Now the three new families of AMS fonts :
%
%   AMS symbol-a
%
\textfont\msafam=\ninemsa \scriptfont\msafam=\sixmsa
\scriptscriptfont\msafam=\fivemsa \def\msa{\fam\msafam\ninemsa}%         
%
%   AMS symbol-b
%
\textfont\msbfam=\ninemsb \scriptfont\msbfam=\sixmsb
\scriptscriptfont\msbfam=\fivemsb \def\msb{\fam\msbfam\ninemsb}%         
%
%   Euler-Fraktal font
%
\textfont\eufmfam=\nineeufm  \scriptfont\eufmfam=\sixeufm
\scriptscriptfont\eufmfam=\fiveeufm \def\frak{\fam\eufmfam\nineeufm}%
%
%%%  Add further fonts families here if you need them.
%
%    Reset \normalbaselineskip and \strubox :
%
\normalbaselineskip=11pt%
\setbox\strutbox=\hbox{\vrule height8pt depth3pt width0pt}%
%
%    Set \normalbaselines and \rm (roman) as defaults :
%
\normalbaselines\rm
%
%    Reset some of the basic vertical skips:
%
\stdskip=4pt plus2pt minus2pt    
\medskipamount=\stdskip          
\parskip=\stdskip                
\abovedisplayskip=\stdskip       
\belowdisplayskip=\stdskip       
\def\ppar{\par\goodbreak\vskip 6pt plus 3pt minus 3pt}%     
%
%   And finally reset the size of section heads (see below):
%
\def\section##1{\global\advance\sectionnumber by 1
\vskip-\lastskip\penalty-800\vskip 20pt plus10pt minus5pt 
\egroup{\bf\number\sectionnumber\quad##1}\bgroup\small         
\vskip 6pt plus3pt minus3pt
\nobreak\resultnumber=1}%      % Reset resultnumber at start of section
}    %%%   End of  \small  macro      
%
%   Two useful abbreviations to keep track of \small material:
\def\beginsmall{\bgroup\small}
\let\endsmall\egroup
%
%
%    The \large  macro is similar (comments abbreviated):
%
%
\def\large{%
\textfont0=\thirtnrm \scriptfont0=\ninerm \scriptscriptfont0=\sevenrm
\def\rm{\fam0\thirtnrm}%
\textfont1=\thirtni \scriptfont1=\ninei \scriptscriptfont1=\seveni
\textfont2=\thirtnsy \scriptfont2=\ninesy \scriptscriptfont2=\sevensy
\textfont3=\thirtnex \scriptfont3=\thirtnex \scriptscriptfont3=\thirtnex
\textfont\bffam=\thirtnbf \scriptfont\bffam=\ninebf
\scriptscriptfont\bffam=\sevenbf \def\bf{\fam\bffam\thirtnbf}%
\textfont\itfam=\thirtnit \def\it{\fam\itfam\thirtnit}%
\textfont\slfam=\thirtnsl \def\sl{\fam\slfam\thirtnsl}%
\textfont\ttfam=\thirtntt \def\tt{\fam\ttfam\thirtntt}%
%   AMS symbol-a  :
\textfont\msafam=\thirtnmsa \scriptfont\msafam=\ninemsa
\scriptscriptfont\msafam=\sevenmsa \def\msa{\fam\msafam\thirtnmsa}%         
%   AMS symbol-b  :
\textfont\msbfam=\thirtnmsb \scriptfont\msbfam=\ninemsb
\scriptscriptfont\msbfam=\sevenmsb \def\msb{\fam\msbfam\thirtnmsb}%         
%   Euler-Fraktal font :
\textfont\eufmfam=\thirtneufm  \scriptfont\eufmfam=\nineeufm
\scriptscriptfont\eufmfam=\seveneufm \def\frak{\fam\eufmfam\teneufm}%
%%%% Add further fonts families here if you need them.
%   Reset \normalbaselineskip and \strubox and initialise :
\normalbaselineskip=16pt%
\setbox\strutbox=\hbox{\vrule height11.5pt depth4.5pt width0pt}%
\normalbaselines\rm}%
\let\Large\large   %  for compatibility with latex
%
%   The next two lines define commonly used switches for
%   blackboard bold (\Bbb) and gothic type (\goth).  The   
%   \Bbb  switch is set to work in the same way as in amstex
%   and switches only the next character to blackboard bold.
\def\Bbb#1{{\msb#1}}

%
%   To use the new AMS fonts you can either use the control
%   sequences \msa, \msb (alias \Bbb) and \frak (alias \goth) eg :

   % see the msam font table
%
%   or, more generally, make \mathchardef's (cf Knuth p155) eg :
\mathchardef\plussquare="0\hexa01
\mathchardef\nge="3\hexb0B
\mathchardef\maltesecross="0\hexa7A
\mathchardef\del="0\hexf01
%
%   or you can use the amstex names for all the new symbols by
%   inserting the line  \input amsnames  in your file directly
%   after \input gtmacros. 
%   This presupposes that you have collected a copy of the file
%   amsnames.tex  from the  gt/info/macros  ftp directory.
%
%
%   Finally we need a small capital font (for author(s)) :
%
\font\sc=cmcsc10
%
%%%%%%%%%%%%%%%%%       END OF FONT MACROS     %%%%%%%%%%%%%
%
%
%                 Knuth's \square macro :
%
\def\sqr#1#2{{\vcenter{\vbox{\hrule  height.#2truept
	\hbox{\vrule width.#2truept height#1truept 
	\kern#1truept \vrule width.#2truept}
	\hrule height.#2truept}}}}
\def\sq{\sqr55}    %   A small square for end-of-proofs. 
%                  %   (Define other size squares by varing the
%                  %   the two numbers.)
%
%
%      Style macros for section heads, theorem statements etc :
%   
%
\newcount\sectionnumber            %%%  Allocate registers to take
\newcount\resultnumber             %%%  section and result numbers.
\sectionnumber=0\resultnumber=1    %%%  Set these registers to 0 and 1
%
%   The \section macro produces a \large bold faced section heading
%   numbered to the left.  Pagebreaks are encouraged before the
%   start of the section and discouraged directly after the heading.
%   Typical use  \section{First steps}  with typical result :
%
%    1  First Steps     (set bold and \large)
%
\def\section#1{\global\advance\sectionnumber by 1
\xdef\nextkey{\number\sectionnumber}%      (used by the \key macro)
\vskip-\lastskip\penalty-800\vskip 20pt plus10pt minus5pt 
{\large\bf\number\sectionnumber\quad#1}         
\vskip 8pt plus4pt minus4pt
\nobreak\resultnumber=1}      % Reset resultnumber at start of section
%
%
%
%   Next a macro to set subheadings (like the  \section  macro
%   but without the number, with less space and set in standard size).
%
%   Typical use :  \sh{Example formats}
%
         
%
%   The \proc ... \endproc macros ("proclaim") are for setting theorems, 
%   lemmas, conjectures etc with automatic numbering.  Typical use :    
%  
%    \proc{Theorem}Every lemon is yellow.\endproc
%
%   Typical result :
%     
%    Theorem 3.4  Every lemon is yellow.   

%   (with Theorem 3.4 set bold and a \stdspace of space before the 
%   statement set in slanted type).
%
\def\proc#1{\xdef\nextkey{\number\sectionnumber.\number\resultnumber}%
\vskip-\lastskip\ppar\bf%
\noindent#1\ \number\sectionnumber.\number\resultnumber
\stdspace\sl\global\advance\resultnumber by 1\ignorespaces}
\def\endproc{\rm\ppar} 
%
%  The \prf ... \endprf macros are for setting proofs.  The code
%  for \prf includes the code for \endproc, so there is no need to
%  type \endproc if the theorem is followed immediatedly by a proof.
%
\def\prf{\vskip-\lastskip\ppar\noindent{\bf Proof}%
\stdspace\rm}                            %  For start of proofs  
   %  For end (or absence) of proofs
\def\endprf{\unskip\stdspace\hbox{}%     %  For end of proof (with
\hfill$\sq$\par\medskip}                 %  extra vertical space)  
        %  For start of proof with alternative name
              %  \endproof is an alias for \endprf
%
%   Typical uses :    
%  
%    \proc{Theorem}Every lemon is yellow. \qed\endproc
%
%    \proc{Theorem}Every lemon is yellow.
%    \prf Use your eyes. \endprf
%
%    \proc{Theorem}Every lemon is yellow.
%    \proof{Proof of theorem} Use your eyes. \endprf
%
%   The next macro is a variant of the \proc macro.  It has
%   exactly the same result except that it omits the number.
%
%   Typical use :  
%    
%    \proclaim{Conjecture}Some oranges are yellow.\endproc
%
\def\proclaim#1{\vskip-\lastskip\ppar\bf%
\noindent#1\stdspace\sl\ignorespaces} 

%
%   The next macro is a further variant for remarks, definitions etc.   
%   It omits the number and does not switch on slanted type.  
%  
%   Typical use :
%
%    \rk{Remark}Some lemons are thick-skinned.\endrk
%
\def\rk#1{\vskip-\lastskip\ppar{\bf #1}\stdspace\ignorespaces}                

%
%   The next macro is for numbering equations etc, \label  produces the 
%   correct number  x.y  and advances the resultnumber register
%
%   Typical use :
%
%     $$fx=7\eqno{\bf\label}$$
%
%   result :
%
%                           fx = 7                           3.5
%
\def\label{\xdef\nextkey{\number\sectionnumber.\number\resultnumber}%
\number\sectionnumber.\number\resultnumber
\global\advance\resultnumber by 1}
%
%
%
%   The next macros are to automate external references.  To use them 
%   type \reflist ..... \endreflist near the beginning of your paper, 
%   where  .... is the list of references in alphabetical order 
%   and in  the form  \key{KEY}  reference    where "KEY" is a 
%   string of characters which reminds you of the reference.   
%   Separate  references with a blank line or a \par.   Eg 
%
%     \reflist
%
%     ..... more references ....
%
%     \key{Kn-84} {\bf D Knuth}, {\it The TeXbook}, Addison--Wesley (1984)
%
%     ..... more references ....
%
%     \endreflist
%
%   Then type  \references  where you wish the references to be printed
%   (normally near the end of the paper).  To refer to Knuth type
%   for example    see Knuth [\ref{Kn-84}, page 320]   and the correct
%   numerical reference will be printed.  Edit the \references macro
%   to change the formatting (if desired).
%   There is an alternative \refkey for \key, provided your KEY contains
%   only letters.  The syntax is:
%
%     \reflist
%
%     ..... more references ....
%
%     \refkey\Knuth  {\bf D Knuth}, {\it The TeXbook}, Addison--Wesley (1984)
%
%     ..... more references ....
%
%     \endreflist
%
%   \key{Knuth}  has exactly the same maening as \refkey\Knuth and you
%   can mix the two syntaxes if you want.  But \refkey\Kn-84
%   would not work.  It would set Kn as the KEY and -84 would get printed!
%
\newcount\refnumber              %  Register for reference numbers
\refnumber=1                     %  set initially to 1.
\long\def\reflist#1\endreflist{%
\long\def\thereflist{#1}{\def\refkey##1##2\par{\xdef##1{\number\refnumber}%
\global\advance\refnumber by 1}%
\def\key##1##2\par{\expandafter\xdef%
\csname##1\endcsname{\number\refnumber}%
\global\advance\refnumber by 1}#1\par}}
\long\def\references{%
\penalty-800\vskip-\lastskip\vskip 15pt plus10pt minus5pt 
{\large\bf References}\ppar %`References' is set \large bold with space around.
{\leftskip=25pt\frenchspacing    % The list of references is set 
\small\parskip=3pt plus2pt       % \small  with small spaces between,
\def\refkey##1##2\par{\noindent  % numbers in [,]'s and set just to the
\llap{[##1]\stdspace}\ignorespaces##2\par}         % left of a 25pt margin.
\def\key##1##2\par{\noindent  
\llap{[\ref{##1}]\stdspace}\ignorespaces##2\par}  
\def\,{\thinspace}\thereflist\par}}
%
%   Next a footnote macro (with automatic numbering) which sets the
%   footnote  \small.
%
%   Typical use :
%         ..... are yellow.\fnote{By yellow here we mean Britsh
%    Standard colour BS3320.} 
%
\newcount\footnotenumber         % Register for footnote number
\footnotenumber=1                % set initially to 1
\def\fnote#1{\xdef\nextkey{\number\footnotenumber}%
{\small\ifnum\footnotenumber>9\parindent=14pt%
\else\parindent=10pt\fi\footnote{$^{\number\footnotenumber}$}%
{\hglue-5pt#1}\global\advance\footnotenumber by 1}}
%
%
%   Next macros for handling figures with automatic numbering (using 
%   TeX's \midinsert to float the figure to a suitable place).
%   
%   The \figure ... \endfigure macro centres the figure and adds
%   an automatically numbered label  Figure XX  after it.
%
%   If you have a caption, then type \caption{caption text} 
%   somewhere between \figure and \endfigure.  The macro
%   will then add  Figure XX: caption text  after the figure.
%
%   If you want an unnumbered or uncentred figure, then use TeX's raw 
%       \midinsert Figure instructions \endinsert  
%   and if you want a numbered figure label in the same style then
%   use \caption{caption text} outside of  \figure ... \endfigure.
%
%   If you need just the label Figure XX  outside of  \figure ... \endfigure
%   then type  \figurelabel .
%
\newcount\figurenumber          % register for figure number
\figurenumber=1                 % set initially to 1
\def\caption#1{\xdef\nextkey{\number\figurenumber}%
\cl{\small Figure \number\figurenumber: #1}%
\global\advance\figurenumber by 1}
\def\figurelabel{\xdef\nextkey{\number\figurenumber}%
\cl{\small Figure \number\figurenumber}%
\global\advance\figurenumber by 1}
\long\def\figure#1\endfigure{{\xdef\nextkey{\number\figurenumber}%
\let\captiontext\relax\def\caption##1{\xdef\captiontext{##1}}%
\midinsert\cl{\ignorespaces#1\unskip\unskip\unskip\unskip}\vglue6pt\cl{\small 
Figure \number\figurenumber\ifx\captiontext\relax\else: \captiontext
\fi}\endinsert\global\advance\figurenumber by 1}}
%
%   Macros for self-correcting internal references.
%
%   There are two macros  \key{KEY}  and  \ref{KEY} .
%
%   The \key macro sets up KEY as a key for whatever number is 
%   being referenced and the \ref macro converts the KEY into 
%   that number.  Type \key after a  \section or \proc or 
%   \label or \fnote or \figure or \caption or \figurelabel .
%
%   Example:
%
%       \section{Introduction}\key{intro}
%       \proc{Theorem}\key{MainTh}Lemons are yelloy\endproc
%       Here we follow\fnote{Follow in the sense of Dickens}
%       \key{Dickens-note}the crowd ....  
%
%       In section \ref{intro}
%       we stated theorem \ref{mainTh} and noted (see footnote 
%       \ref{Dickens-note}) ...
%
\def\nextkey{??}   %  initialise \nextkey (which is reset by all the
%                     numbering macros)
%
\def\key#1{\expandafter\xdef\csname #1\endcsname{\nextkey}}
\def\ref#1{\expandafter\ifx\csname #1\endcsname\relax
\immediate\write16{Reference {#1} undefined}??\else
\csname #1\endcsname\fi}
%
%   Note:  If the KEY contains only letters then \KEY has exactly the
%   same meaning as \ref{KEY} so in the example you could have:
%
%       In section \intro\ we ....
%
%   The \key will work at any time after the macro which sets the
%   number, provided no other macro which sets a number has been used. 
%
%   Macros for forward references:
%              =======
%   The \key \ref macros ONLY work for backwards references.  If you  
%   want to use forwards references, then type \useforwardrefs  near
%   the beginning of your file.  The KEY's are then stored in an
%   auxiliary  .ref  file and you then suffer the same disadvantage as
%   when using LaTeX that you must TeX the file twice to get
%   the references correct.
%
%   To use a forward ref type \ref{KEY}.  (You can type the
%   alternative  \KEY  but you'll get an error on first TeX'ing 
%   if the \KEY is not yet defined.) 
%
%   The macro also allows external references to be listed at the end 
%   of the file (if you wish to).  (Indeed they can be typed anywhere
%   before the \references command.)  You can combine the reference list
%   and the \references command by typing the references (using the
%   same syntax as before) between the commands \biblio and \endbiblio 
%   (don't type \references or they'll be printed twice).
%
\newread\gtinfile
\newwrite\gtreffile
\def\useforwardrefs{
\openin\gtinfile\jobname.ref
\ifeof\gtinfile
\closein\gtinfile
\immediate\write16{No file \jobname.ref}
\else
\closein\gtinfile
\input \jobname.ref
\fi
\immediate\openout\gtreffile \jobname.ref
%
%   Adapt \key :
%
\def\key##1{{\def\\{\noexpand}%
\expandafter\xdef\csname ##1\endcsname{\nextkey}%
\immediate\write\gtreffile{\\\expandafter\\\def\\\csname ##1\\\endcsname%
{\nextkey}}}}
%
%  Adapt macros for external references:  
%
\long\def\reflist##1\endreflist{%
\long\def\thereflist{##1}{\def\refkey####1####2\par{\xdef####1{%
\number\refnumber}{\def\\{\noexpand}\immediate\write\gtreffile
{\\\def\\####1{\number\refnumber}}}\global\advance\refnumber by 1}%
\def\key####1####2\par{\expandafter\xdef%
\csname####1\endcsname{\number\refnumber}%
{\def\\{\noexpand}\immediate\write\gtreffile
{\\\expandafter\\\def\\\csname ####1\\\endcsname{\number\refnumber}}}
\global\advance\refnumber by 1}##1\par}}
\long\def\biblio##1\endbiblio{\reflist##1\endreflist\references}%
%
%  Adapt obselete key macros (\numkey, \seckey and \figkey):
%
\def\numkey##1{{\def\\{\noexpand}%
\xdef##1{\number\sectionnumber.\number\resultnumber}
\immediate\write\gtreffile{\\\def\\##1%
{\number\sectionnumber.\number\resultnumber}}}}
\def\seckey##1{{\def\\{\noexpand}\xdef##1{\number\sectionnumber}
\immediate\write\gtreffile{\\\def\\##1{\number\sectionnumber}}}}
\def\figkey##1{\xdef##1{\number\figurenumber}%
{\def\\{\noexpand}\immediate\write\gtreffile%
{\\\def\\##1{\number\figurenumber}}}
\number\figurenumber\global\advance\figurenumber by 1}
}   %  end of \useforwardrefs
%
%
%   The next five macros are obselete and have been superseeded by
%   the general \key macro above.  They are included merely to 
%   maintain backward compatibility for the package:
%
%
\def\figkey#1{\xdef#1{\number\figurenumber}%
\number\figurenumber\global\advance\figurenumber by 1}
\def\fig#1#2\endfig{%
\midinsert\cl{#2}\vglue6pt\cl{\small Figure #1}\endinsert}
\def\newfig{\number\figurenumber\global\advance\figurenumber by 1}
\def\numkey#1{\xdef#1{\number\sectionnumber.\number\resultnumber}}
\def\seckey#1{\xdef#1{\number\sectionnumber}}
%
%   End of obselete macros.
%
%
%   The next macro is a version of the verbatim macro given by Knuth.
%
%   This macro produces a "verbatim" printout of
%   any ASCII string which does not contain the symbol "
%   (TeX files do not usually contain " 's).
%   More precisely, everything between consecutive pairs
%   of " 's is printed verbatim in the typewriter font cmtt.
%   For an explanation of how the macro works, see Knuth pp 420-1.
%
%   There are two switches: \verb (which switches the macro on)
%   and \brev which switches the macro off (the default).  When
%   the macro is switched off the symbol " has its usual 
%   meaning for TeX.  To use the macro, type \verb before use
%   and the use " to switch verbatim on and off.  Be careful
%   not to use " for any other purpose.  There is no need to
%   switch the macro off again unless you need to use " for
%   some other purpose (eg making  \mathchardef 's).  Note 
%   that the macro MUST BE OFF before inputting  amsnames.tex .
%
%   Whether the macro is on or off you can always use the
%   control sequence \dq (double quote) for " e.g.
%   \mathchardef\sum=\dq1350  is perfectly valid.
%   The control sequence \ttq is an abbreviation for
%   {\tt\dq}.  Thus "\ttq" will produce " (in cmtt)
%   inside a verbatim quote.
%
%
   %  define a code for " so it can be used when \verb is on
  %  code for " in cmtt
%
\def\verb{\catcode`\"=\active}       %  The main
\def\brev{\catcode`\"=12}            %  switches.
\brev                                %  Prime switches and
\verb                                %  switch on.
{\obeyspaces\gdef {\ }}              
{\catcode`\`=\active\gdef`{\relax\lq}}
\def"{%
\begingroup\baselineskip=12pt\def\par{\leavevmode\endgraf}%
\tt\obeylines\obeyspaces\parskip=0pt\parindent=0pt%
\catcode`\$=12\catcode`\&=12\catcode`\^=12\catcode`\#=12%
\catcode`\_=12\catcode`\~=12%
\catcode`\{=12\catcode`\}=12\catcode`\%=12\catcode`\\=12%
\catcode`\`=\active\let"\endgroup}
\brev      %   Finally switch the macro off (for safety)
%
%   Macros for itemised lists.   Typical use :
%    
%    \items
%    \item{(i)}Colours must be defined.
%    \item{(ii)}Colour cards may not be cited.
%    \enditems
%
%   Result :
%
%    (i)  Colours must be defined. 
%   (ii)  Colour cards may not be cited.
%
%
\def\items{\par\leftskip = 25pt}           % Start of itemised list         
\def\enditems{\par\leftskip = 0pt}         % end of itemised list   
\def\item#1{\par\leavevmode\llap{#1\stdspace}%
\ignorespaces}                             % labelled item
\def\itemb{\item{$\bullet$}}               % bulleted item.
%
%   The \quote ... \endquote macros are for typesetting quotations :
%

%
%   A few useful abbreviations :
%
    %  Colon with correct spacing for maps.
\def\np{\vfil\eject}         %  Forced page break (new page).
\def\nl{\hfil\break}         %  New line.
\def\cl{\centerline}         %  Centerline
\def\gt{{\mathsurround=0pt\it $\cal G\mskip-2mu$eometry \&\ 
$\cal T\!\!$opology}}        %  The journal title in recommended style
    %  for monographs
\def\agt{{\mathsurround=0pt\it$\cal A\mskip-.7mu$lgebraic \&\ 
$\cal G\mskip-2mu$eometric $\cal T\!\!$opology}}  % AGT
%
%    Finally some macros for automatic title page or header generation.
%    To use them type your header information using the following  
%    example as a guide :
%
%    Note that \\ is used as standard separator (for lines in \title and
%    \address, between authors and between email addresses or URL's)
%    and that \email, \url and \secondaddress are optional.
%

% Example:  \title{A short spoof paper\\with a two-line title}
% =======   \authors{Albert Einstein\\Leonardo da Vinci}
%           \address{IAS\\Princeton}\secondaddress{Renaissance\\Venice}
%           \email{ae@ias.princeton.edu\\ldv@ren.ven.hist}
%           \abstract 
%           A short spoof paper with a very short abstract.
%           \endabstract 
%           \primaryclass{00-01, 00-02}\secondaryclass{68-00, 68-01}
%           \keywords{Short, spoof, paper}
%           \maketitlepage
%
%
%    The title page or header will then be generated automatically.
%
%
%    Define the various ingredients of the title page:
%
\def\title#1{\def\thetitle{#1}}
\def\shorttitle#1{\def\theshorttitle{#1}}
\def\author#1{\edef\previousauthors{\theauthors}
 \ifx\theauthors\relax\def\theauthors{#1}\else
 \def\theauthors{\previousauthors\par#1}\fi}

\let\authors\author        % aliases
\def\address#1{\edef\previousaddresses{\theaddress}
 \ifx\theaddress\relax\def\theaddress{#1}\else
 \def\theaddress{\previousaddresses\par\vskip 2pt\par#1}\fi}
                             % alias
\def\secondaddress#1{\edef\previousaddresses{\theaddress}
 \ifx\theaddress\relax\def\theaddress{#1}\else
 \def\theaddress{\previousaddresses\par{\rm and}\par#1}\fi}   

\def\email#1{\edef\previousemails{\theemail}
 \ifx\theemail\relax\def\theemail{#1}\else
 \def\theemail{\previousemails\hskip 0.75em\relax#1}\fi}
  % aliases
\def\secondemail#1{\edef\previousemails{\theemail}
 \ifx\theemail\relax\def\theemail{#1}\else
 \def\theemail{\previousemails\hskip 0.75em{\rm and}\hskip 0.75em
 \relax#1}\fi}
\def\url#1{\edef\previousurls{\theurl}
 \ifx\theurl\relax\def\theurl{#1}\else
 \def\theurl{\previousurls\hskip 0.75em\relax#1}\fi}
      % aliases
\def\secondurl#1{\edef\previousurls{\theurl}
 \ifx\theurl\relax\def\theurl{#1}\else
 \def\theurl{\previousurls\hskip 0.75em{\rm and}\hskip 0.75em
 \relax#1}\fi}
\long\def\abstract#1\endabstract{\long\def\theabstract{#1}}
\def\primaryclass#1{\def\theprimaryclass{#1}}
                        % alias
\def\secondaryclass#1{\def\thesecondaryclass{#1}}
\def\keywords#1{\def\thekeywords{#1}}
%
%  Set \\ to \par and title page items to \relax to initialise macros :
%
\let\\\par\let\thetitle\relax\let\theshorttitle\relax
\let\theauthors\relax\let\theshortauthors\relax
\let\theaddress\relax\let\theshortaddress\relax
\let\theemail\relax\let\theurl\relax
\let\theabstract\relax\let\theprimaryclass\relax
\let\thesecondaryclass\relax\let\thekeywords\relax
%
%
%
%   Basic title page layout (edit this macro if you
%   wish to adjust the title page layout) :
%
\long\def\maketitlepage{    % start of definition of \maketitlepage

\vglue 0.2truein   % top margin

% title :
%
{\parskip=0pt\leftskip 0pt plus 1fil\def\\{\par\smallskip}{\large
\bf\thetitle}\par\medskip}   

\vglue 0.15truein 

% authors :
%
{\parskip=0pt\leftskip 0pt plus 1fil\def\\{\par}{\sc\theauthors}
\par\medskip}%
 
\vglue 0.1truein 

% address(es) email's and URL's (with switches to detect whether the
% optional items have been used) :
%
{\small\parskip=0pt
{\leftskip 0pt plus 1fil\def\\{\par}{\sl\theaddress}\par}
\ifx\theemail\relax\else  % email address?
\vglue 5pt \def\\{\stdspace{\rm and}\stdspace} 
\cl{Email:\stdspace\tt\theemail}\fi
\ifx\theurl\relax\else    % URL given?
\vglue 5pt \def\\{\stdspace{\rm and}\stdspace} 
\cl{URL:\stdspace\tt\theurl}\fi\par}

\vglue 7pt 

{\bf Abstract}

\vglue 5pt

\theabstract

\vglue 7pt 

{\bf AMS Classification numbers}\quad Primary:\quad \theprimaryclass\par

Secondary:\quad \thesecondaryclass

\vglue 5pt 

{\bf Keywords:}\quad \thekeywords

\np  % page break at the end of the title page

}    % end of definition of \maketitlepage
%
%    % \makeshorttitle (for general preprints) doesn't take a new page
%
\long\def\makeshorttitle{    % start of definition of \makeshorttitle

%\vglue 0.2truein   % top margin

% title :
%
{\parskip=0pt\leftskip 0pt plus 1fil\def\\{\par\smallskip}{\large
\bf\thetitle}\par\medskip}   

\vglue 0.05truein 

% authors :
%
{\parskip=0pt\leftskip 0pt plus 1fil\def\\{\par}{\sc\theauthors}
\par\medskip}%
 
\vglue 0.03truein 

% address(es) email's and URL's (with switches to detect whether the
% optional items have been used) :
%
{\small\parskip=0pt
{\leftskip 0pt plus 1fil\def\\{\par}{\sl\ifx\theshortaddress\relax
\theaddress\else\theshortaddress\fi}\par}
\ifx\theemail\relax\else  % email address?
\vglue 5pt \def\\{\stdspace{\rm and}\stdspace} 
\cl{Email:\stdspace\tt\theemail}\fi
\ifx\theurl\relax\else    % URL given?
\vglue 5pt \def\\{\stdspace{\rm and}\stdspace} 
\cl{URL:\stdspace\tt\theurl}\fi\par}

\vglue 10pt 

% abstract and classification numbers (with switches):

{\small\leftskip 25pt\rightskip 25pt{\bf Abstract}\stdspace\theabstract

{\bf AMS Classification}\stdspace\theprimaryclass
\ifx\thesecondaryclass\relax\else; \thesecondaryclass\fi\par
{\bf Keywords}\stdspace \thekeywords\par}
\vglue 7pt
}    % end of definition of \makeshorttitle
\let\maketitle\makeshorttitle        %% alias
%
%    %%%% \makeagttitle (for AGT) similar to \makeshorttitle but
%         with addresses omitted (they go at the end)
%
%%%% publication info and test defaults:

\def\volumenumber#1{\def\thevolumenumber{#1}}
\def\volumeyear#1{\def\thevolumeyear{#1}}
\def\pagenumbers#1#2{\def\startpage{#1}\def\finishpage{#2}}
\def\published#1{\def\publishdate{#1}}
\def\received#1{\def\receiveddate{#1}}
\def\revised#1{\def\reviseddate{#1}}
\let\reviseddate\relax
%% Defaults for authors to use to check layout
\volumenumber{X}
\volumeyear{20XX}
\pagenumbers{1}{XXX}
\published{XX Xxxember 20XX}

\long\def\makeagttitle{   %%% start of definition of \makeagttitle
\agt\hfill      %   Journal title (top left) 
%   logo placeholder (top right)
\hbox to 60truept{\vbox to 0pt{\vglue -14truept{\bf [Logo here]}\vss}\hss}
\break
{\small Volume \thevolumenumber\ (\thevolumeyear)
\startpage--\finishpage\nl
Published: \publishdate}

\vglue .2truein

% title
{\parskip=0pt\leftskip 0pt plus 1fil\def\\{\par\smallskip}{\large
\bf\thetitle}\par\medskip}   
\vglue 0.05truein 

% authors :
%
{\parskip=0pt\leftskip 0pt plus 1fil\def\\{\par}{\sc\theauthors}
\par\medskip}%
 
\vglue 0.03truein 

%  abstract and classification numbers:

{\small\leftskip 25truept\rightskip 25truept{\bf Abstract}\stdspace\theabstract

{\bf AMS Classification}\stdspace\theprimaryclass
\ifx\thesecondaryclass\relax\else; \thesecondaryclass\fi\par
{\bf Keywords}\stdspace \thekeywords\par}\vglue 7truept

}   %%%% end of definition of \makeagttitle

%%%%% Macro to typeset addresses (typically at the end of the paper)

\def\Addresses{\bigskip
{\small \parskip 0pt \leftskip 0pt \rightskip 0pt plus 1fil \def\\{\par}
\sl\theaddress\par\medskip \rm Email:\stdspace\tt\theemail\par
\ifx\theurl\relax\else\smallskip \rm URL:\stdspace\tt\theurl\par\fi}}

\def\agtart{%   Full mock-up of AGT article style (for authors to test with)
%  get print centerpage:
\hoffset 14truemm
\voffset 31truemm
\font\phead=cmsl9 scaled 950
\font\pnum=cmbx10 scaled 913
\font\pfoot=cmsl9 scaled 950
%  headline and footline
\headline{\vbox to 0pt{\vskip -4.5mm\line{\small\phead\ifnum
\count0=\startpage ISSN numbers are printed here
\hfill {\pnum\folio}\else\ifodd\count0\def\\{ }% 
\ifx\theshorttitle\relax\thetitle\else\theshorttitle\fi\hfill{\pnum\folio}
\else\def\\{ and }{\pnum\folio}\hfill\ifx\theshortauthors\relax\theauthors
\else\theshortauthors\fi\fi\fi}\vss}}
\footline{\vbox to 0pt{\vglue 0mm\line{\small\pfoot\ifnum\count0=\startpage
Copyright declaration is printed here\hfill\else
\agt, Volume \thevolumenumber\ (\thevolumeyear)\hfill\fi}\vss}}
%  force \agttitle
\let\maketitle\makeagttitle\let\makeshorttitle\makeagttitle}

\input epsf
\chardef\newinsCatAt\the\catcode `\@
\catcode `\@=11
%
%%%%%%%%%%%% Corrected insert macros for plain.tex %%%%%%%%%%%%%%%%%
%
%  New skipamounts:
%
\newskip\insertskipamount\newskip\inserthardskipamount
\insertskipamount 12pt plus2pt     % Redefined (CPR)    % default 6pt plus 2pt
\inserthardskipamount 4pt          % to suit GT style   % default 6pt
\def\insertskip{\vskip\insertskipamount}
%
%  Save and restore \lastskip:
%
\newskip\LastSkip
\def\SaveLastSkip{\LastSkip\lastskip}
\def\RestoreLastSkip{\nobreak\vskip-\LastSkip\vskip\LastSkip}
%
%  Larry Siebenmann's test for split topinserts:
%
\newcount\SplitTest%        will be set to -1 if a topinsert has split
\def\SetSplitTest{\SplitTest\insertpenalties
  \insert\topins{\floatingpenalty1}%
  \advance\SplitTest-\insertpenalties}
%
%  From here on we modify definitions in plain.tex.
%
% Redefine \midinsert to convert to \topinsert if a topinsert has been
% split, to prevent midinserts getting out of order (cf. TeXbook Exercise
% 15.5). As in plain.tex, a \midinsert still converts to a \topinsert
% (which then splits) if the insert is too big for current page.
%   Was:    \def\midinsert{\@midtrue\@ins}
\def\midinsert{\par
 \SaveLastSkip\penalty-150\SetSplitTest\RestoreLastSkip
 \ifnum\SplitTest=-1
  \@midfalse\p@gefalse\else\@midtrue\fi\@ins}
% Redefine \@ins to add \inserthardskipamount of glue above.
%   Was:  \def\@ins{\par\begingroup\setbox\z@\vbox\bgroup}
\def\@ins{\par\begingroup\setbox\z@\vbox\bgroup%
  \vglue\inserthardskipamount}
% Changes to \endinsert of plain.tex 3.0:
% - Use \insertskipamount instead of \bigskipamount throughout.
% - Use larger of previous skip and insertskip before middle insert.
% - Add \nointerlineskip to avoid unwanted extra 1pt skip.
% - Save and restore lastskip when an insert floats.
\def\endinsert{\egroup % finish the \vbox
  \if@mid \dimen@\ht\z@ \advance\dimen@\dp\z@
    \advance\dimen@\insertskipamount%            was 12pt (wn)
    \advance\dimen@\pagetotal\advance\dimen@-\pageshrink
    \ifdim\dimen@>\pagegoal\@midfalse\p@gefalse\fi\fi
% Next 3 lines replace:  \if@mid \bigskip\box\z@\bigbreak (wn)
  \if@mid%
    \ifdim\lastskip<\insertskipamount\removelastskip\insertskip\fi
    \nointerlineskip\box\z@\penalty-200\insertskip
  \else%
    \SaveLastSkip%                                  added (wn)
    \insert\topins{\penalty100 % floating insertion
    \splittopskip\z@skip
    \splitmaxdepth\maxdimen \floatingpenalty\z@
    \ifp@ge \dimen@\dp\z@
    \vbox to\vsize{\unvbox\z@\kern-\dimen@}% depth is zero
    \else \box\z@\nobreak\insertskip\fi}% was \bigskip\fi (wn)
    \RestoreLastSkip%                               added (wn)
   \fi\endgroup}
%%%%%%%%%%%%%%%%%% Done correcting insert macros %%%%%%%%%%%%%%%%%%%
%
\catcode `\@=\newinsCatAt

%%%
%%%  This version of  gtoutput.tex  is intended to finish formatting
%%%  papers published in Geometry & Topology and stored in the
%%%  arXiv.   All versions of  gtoutput.tex  are copyright 
%%%  GT Publications and are to be used _only_ for formatting
%%%  the officially published version of G&T papers.
%%%
%%%
%%%                                             Colin Rourke  14.9.2000
%%%
%%%  To create header file  head.xxx  comment out the first \endinput

%  test for latex or plain tex
\def\ifplaintex{\expandafter\ifx\csname documentclass\endcsname\relax}

%  get print centerpage:

\ifplaintex 
\hoffset 14truemm
\voffset 31truemm
\else
\headsep 23pt
\footskip 35pt
\hoffset -4truemm
\voffset 12.5truemm
\fi

%  load pictex if not already loaded :
\expandafter\ifx\csname beginpicture\endcsname\relax
\expandafter\ifx\csname documentclass\endcsname\relax
\input pictex \else\font\fiverm=cmr5
\input prepictex \input pictex \input postpictex \fi\fi

\def\gt{{\mathsurround=0pt\it $\cal G\mskip-2mu$eometry \&\ 
$\cal T\!\!$opology}}        %  journal title in recommended style

\def\gtp{{\mathsurround=0pt\it $\cal G\mskip-2mu$eometry \&\ 
$\cal T\!\!$opology $\cal P\!$ublications}}  % GT publications

%  define the various new ingredients of the title page 

\def\volumenumber#1{\def\thevolumenumber{#1}}
\def\papernumber#1{\def\thepapernumber{#1}}
\def\volumeyear#1{\def\thevolumeyear{#1}}

\def\pagenumbers#1#2{\def\startpage{#1}\def\finishpage{#2}}
\def\published#1{\def\publishdate{#1}}
\def\proposed#1{\def\theproposer{#1}}
\def\seconded#1{\def\theseconders{#1}}
\def\received#1{\def\receiveddate{#1}}
\def\revised#1{\def\reviseddate{#1}}
\def\accepted#1{\def\accepteddate{#1}}
\def\asciititle#1{\def\theasciititle{#1}}

\def\asciiaddress#1{\def\theasciiaddress{#1}}

\long\def\asciiabstract#1{\long\def\theasciiabstract{#1}}

\def\shorttitle#1{\def\theshorttitle{#1}}

%  initialise

\let\\\par
\let\thevolumenumber\relax\let\thepapernumber\relax
\let\thevolumeyear\relax\let\thesamplenumber\relax\let\startpage\relax
\let\finishpage\relax\let\publishdate\relax\let\receiveddate\relax
\let\reviseddate\relax\let\accepteddate\relax\let\theasciititle\relax
\let\theasciiauthors\relax\let\theasciiaddress\relax
\let\theasciiabstract\relax
\let\theasciiemail\relax\let\theshortauthors\relax\let\theshorttitle\relax

\long\def\maketitlep{   % start of definition of \maketitlep

\count0=\startpage

\gt\hfill      %   Journal title (top left) 
%    Logo (top right) :
\beginpicture
\setcoordinatesystem units <0.33truein, 0.33truein> point at 2.2 0.9
\setplotsymbol ({$\cal G$})
\plotsymbolspacing=9truept
\circulararc 315 degrees from 0 1 center at 0 0
\setplotsymbol ({$\cal T$})
\circulararc 315 degrees from 1 -1 center at 1 0
\endpicture
%   end of logo
%
\break
{\small\ifx\thesamplenumber\relax % sample?  
Volume \else Sample
\fi\thevolumenumber\ (\thevolumeyear)
\startpage--\finishpage\nl
Published: \publishdate}
\vglue 0.5truein plus 0.4fil minus 0.1truein

% title
{\parskip=0pt\leftskip 0pt plus 1fil\def\\{\par\smallskip}{\ifplaintex\large
\else\Large\fi\bf\thetitle}\par\medskip}   

\vglue 0pt plus 0.1fil 

% authors
{\parskip=0pt\leftskip 0pt plus 1fil\def\\{\par}{\sc\theauthors}
\par\medskip}

\vglue 0pt plus 0.1fil 

%address(es)
{\small\parskip=0pt\let\newline\\
{\leftskip 0pt plus 1fil\def\\{\par}{\sl\theaddress}\par}
\expandafter\ifx\theemail\relax    % email address?
\relax\else\vglue 5pt plus 0.02fil minus 2pt\def\\{\stdspace{\rm 
and}\stdspace} 
\cl{Email:\stdspace\tt\theemail}\fi
\ifx\theurl\relax                  % URL given?
\relax\else\vglue 5pt plus 0.02fil minus 2pt\def\\{\stdspace{\rm 
and}\stdspace}
\cl{URL:\stdspace\tt\theurl}\fi\par}

\vglue 7pt plus 0.3fil minus 3pt

{\bf Abstract}
\vglue 5pt plus 0.1fil minus 2pt

\theabstract

\vglue 7pt plus 0.3fil minus 3pt

{\bf AMS Classification numbers}\quad Primary:\quad \theprimaryclass

Secondary:\quad \thesecondaryclass

\vglue 5pt plus 0.3fil minus 2pt

{\bf Keywords:}\quad \thekeywords

\vglue 10pt plus 0.5fil minus 5pt

{\small  Proposed: \theproposer\hfill Received: \receiveddate\nl
Seconded: \theseconders\hfill 
\ifx\reviseddate\relax                         % paper revised?
Accepted: \accepteddate                        % no
\else
Revised: \reviseddate                          % yes
\fi}
\eject
}       %  end of definition of \maketitlep

\let\maketitlepage\maketitlep
\let\maketitle\maketitlepage

%%% Headers and footers

\font\phead=cmsl9 scaled 950
\font\lhead=cmsl9 scaled 1050
\font\pnum=cmbx10 scaled 913
\font\lnum=cmbx10 
\font\pfoot=cmsl9 scaled 950
\font\lfoot=cmsl9 scaled 1050
\ifplaintex
\headline{\vbox to 0pt{\vskip -4.5mm\line{\small\phead\ifnum
\count0=\startpage ISSN 1364-0380 (on line)
1465-3060 (printed) \hfill {\pnum\folio}\else\ifodd\count0\def\\{ }% 
\ifx\theshorttitle\relax\thetitle\else\theshorttitle\fi\hfill{\pnum\folio}
\else\def\\{ and }{\pnum\folio}\hfill\ifx\theshortauthors\relax\theauthors
\else\theshortauthors\fi\fi\fi}\vss}}
\footline{\vbox to 0pt{\vglue 0mm\line{\small\pfoot\ifnum\count0=\startpage
\copyright\ \gtp\hfill\else
\gt, Volume \thevolumenumber\ (\thevolumeyear)\hfill\fi}\vss
}}
\else
\makeatletter
\def\@oddhead{{\small\lhead\ifnum\count0=\startpage ISSN 1364-0380 (on line)
1465-3060 (printed) \hfill {\lnum\number\count0}\else\ifodd\count0
\def\\{ }\ifx\theshorttitle\relax \thetitle \else\theshorttitle\fi\hfill
{\lnum\number\count0}\else\def\\{ and }{\lnum\number\count0}
\hfill\ifx\theshortauthors\relax 
\theauthors\else\theshortauthors\fi\fi\fi}}\def\@evenhead{\@oddhead}
\def\@oddfoot{\small\lfoot\ifnum\count0=\startpage\copyright\ \gtp\hfill\else
\gt, Volume \thevolumenumber\ (\thevolumeyear)\hfill\fi}
\def\@evenfoot{\@oddfoot}
\makeatother
\fi

\newwrite\gtoutfile
\long\gdef\makeheadfile{  %%% start of definition of \makeheadfile
{\def\\{, }\def\s{ }
\immediate\openout\gtoutfile head.xxx
\immediate\write\gtoutfile{To: math@arxiv.org}
\immediate\write\gtoutfile{Subject: put or rep NNNNN:pppp}
\immediate\write\gtoutfile{--text follows this line--}
\immediate\write\gtoutfile{Proxy-for: \ifx\theasciiauthors\relax
\theauthors\else\theasciiauthors\fi\s<\ifx\theasciiemail\relax\theemail\else\theasciiemail\fi>}
\immediate\write\gtoutfile{\noexpand\\}
\immediate\write\gtoutfile{Authors: \ifx\theasciiauthors\relax
\theauthors\else\theasciiauthors\fi}
\immediate\write\gtoutfile{Title: \ifx\theasciititle\relax
\thetitle\else\theasciititle\fi}
\immediate\write\gtoutfile{Subj-class: GT or SG or MG etc}
\immediate\write\gtoutfile{MSC-class: \theprimaryclass\ifx\thesecondaryclass\relax\else, \thesecondaryclass\fi}
\immediate\write\gtoutfile{Journal-ref: Geom. Topol. \thevolumenumber
(\thevolumeyear) \startpage-\finishpage}
\immediate\write\gtoutfile{Comments: Published in Geometry and Topology at}
\immediate\write\gtoutfile{\s\s http://www.maths.warwick.ac.uk/gt/GTVol\thevolumenumber/paper\thepapernumber.abs.html}
\immediate\write\gtoutfile{\noexpand\\}
\immediate\write\gtoutfile{}
\ifx\theasciiabstract\relax
\immediate\write\gtoutfile{\theabstract}\else
\immediate\write\gtoutfile{\theasciiabstract}\fi
\immediate\write\gtoutfile{}
\immediate\write\gtoutfile{\noexpand\\}
\immediate\write\gtoutfile{}
\immediate\write\gtoutfile{<uuencoded .tar.gz file here>}
\immediate\write\gtoutfile{}
\immediate\closeout\gtoutfile}}  %%% end of definition of \makeheadfile

\def\maketitlepage{\maketitlep\makeheadfile}
\let\maketitle\maketitlepage

 %%%
%%%  This version of  gtoutput.tex  is intended to finish formatting
%%%  papers published in Geometry & Topology and stored in the
%%%  arXiv.   All versions of  gtoutput.tex  are copyright 
%%%  GT Publications and are to be used _only_ for formatting
%%%  the officially published version of G&T papers.
%%%
%%%
%%%                                             Colin Rourke  14.9.2000
%%%
%%%  To create header file  head.xxx  comment out the first \endinput

%  test for latex or plain tex
\def\ifplaintex{\expandafter\ifx\csname documentclass\endcsname\relax}

%  get print centerpage:

\ifplaintex 
\hoffset 14truemm
\voffset 31truemm
\else
\headsep 23pt
\footskip 35pt
\hoffset -4truemm
\voffset 12.5truemm
\fi

%  load pictex if not already loaded :
\expandafter\ifx\csname beginpicture\endcsname\relax
\expandafter\ifx\csname documentclass\endcsname\relax
\input pictex \else\font\fiverm=cmr5
\input prepictex \input pictex \input postpictex \fi\fi

\def\gt{{\mathsurround=0pt\it $\cal G\mskip-2mu$eometry \&\ 
$\cal T\!\!$opology}}        %  journal title in recommended style

\def\gtp{{\mathsurround=0pt\it $\cal G\mskip-2mu$eometry \&\ 
$\cal T\!\!$opology $\cal P\!$ublications}}  % GT publications

%  define the various new ingredients of the title page 

\def\volumenumber#1{\def\thevolumenumber{#1}}
\def\papernumber#1{\def\thepapernumber{#1}}
\def\volumeyear#1{\def\thevolumeyear{#1}}

\def\pagenumbers#1#2{\def\startpage{#1}\def\finishpage{#2}}
\def\published#1{\def\publishdate{#1}}
\def\proposed#1{\def\theproposer{#1}}
\def\seconded#1{\def\theseconders{#1}}
\def\received#1{\def\receiveddate{#1}}
\def\revised#1{\def\reviseddate{#1}}
\def\accepted#1{\def\accepteddate{#1}}
\def\asciititle#1{\def\theasciititle{#1}}

\def\asciiaddress#1{\def\theasciiaddress{#1}}

\long\def\asciiabstract#1{\long\def\theasciiabstract{#1}}

\def\shorttitle#1{\def\theshorttitle{#1}}

%  initialise

\let\\\par
\let\thevolumenumber\relax\let\thepapernumber\relax
\let\thevolumeyear\relax\let\thesamplenumber\relax\let\startpage\relax
\let\finishpage\relax\let\publishdate\relax\let\receiveddate\relax
\let\reviseddate\relax\let\accepteddate\relax\let\theasciititle\relax
\let\theasciiauthors\relax\let\theasciiaddress\relax
\let\theasciiabstract\relax
\let\theasciiemail\relax\let\theshortauthors\relax\let\theshorttitle\relax

\long\def\maketitlep{   % start of definition of \maketitlep

\count0=\startpage

\gt\hfill      %   Journal title (top left) 
%    Logo (top right) :
\beginpicture
\setcoordinatesystem units <0.33truein, 0.33truein> point at 2.2 0.9
\setplotsymbol ({$\cal G$})
\plotsymbolspacing=9truept
\circulararc 315 degrees from 0 1 center at 0 0
\setplotsymbol ({$\cal T$})
\circulararc 315 degrees from 1 -1 center at 1 0
\endpicture
%   end of logo
%
\break
{\small\ifx\thesamplenumber\relax % sample?  
Volume \else Sample
\fi\thevolumenumber\ (\thevolumeyear)
\startpage--\finishpage\nl
Published: \publishdate}
\vglue 0.5truein plus 0.4fil minus 0.1truein

% title
{\parskip=0pt\leftskip 0pt plus 1fil\def\\{\par\smallskip}{\ifplaintex\large
\else\Large\fi\bf\thetitle}\par\medskip}   

\vglue 0pt plus 0.1fil 

% authors
{\parskip=0pt\leftskip 0pt plus 1fil\def\\{\par}{\sc\theauthors}
\par\medskip}

\vglue 0pt plus 0.1fil 

%address(es)
{\small\parskip=0pt\let\newline\\
{\leftskip 0pt plus 1fil\def\\{\par}{\sl\theaddress}\par}
\expandafter\ifx\theemail\relax    % email address?
\relax\else\vglue 5pt plus 0.02fil minus 2pt\def\\{\stdspace{\rm 
and}\stdspace} 
\cl{Email:\stdspace\tt\theemail}\fi
\ifx\theurl\relax                  % URL given?
\relax\else\vglue 5pt plus 0.02fil minus 2pt\def\\{\stdspace{\rm 
and}\stdspace}
\cl{URL:\stdspace\tt\theurl}\fi\par}

\vglue 7pt plus 0.3fil minus 3pt

{\bf Abstract}
\vglue 5pt plus 0.1fil minus 2pt

\theabstract

\vglue 7pt plus 0.3fil minus 3pt

{\bf AMS Classification numbers}\quad Primary:\quad \theprimaryclass

Secondary:\quad \thesecondaryclass

\vglue 5pt plus 0.3fil minus 2pt

{\bf Keywords:}\quad \thekeywords

\vglue 10pt plus 0.5fil minus 5pt

{\small  Proposed: \theproposer\hfill Received: \receiveddate\nl
Seconded: \theseconders\hfill 
\ifx\reviseddate\relax                         % paper revised?
Accepted: \accepteddate                        % no
\else
Revised: \reviseddate                          % yes
\fi}
\eject
}       %  end of definition of \maketitlep

\let\maketitlepage\maketitlep
\let\maketitle\maketitlepage

%%% Headers and footers

\font\phead=cmsl9 scaled 950
\font\lhead=cmsl9 scaled 1050
\font\pnum=cmbx10 scaled 913
\font\lnum=cmbx10 
\font\pfoot=cmsl9 scaled 950
\font\lfoot=cmsl9 scaled 1050
\ifplaintex
\headline{\vbox to 0pt{\vskip -4.5mm\line{\small\phead\ifnum
\count0=\startpage ISSN 1364-0380 (on line)
1465-3060 (printed) \hfill {\pnum\folio}\else\ifodd\count0\def\\{ }% 
\ifx\theshorttitle\relax\thetitle\else\theshorttitle\fi\hfill{\pnum\folio}
\else\def\\{ and }{\pnum\folio}\hfill\ifx\theshortauthors\relax\theauthors
\else\theshortauthors\fi\fi\fi}\vss}}
\footline{\vbox to 0pt{\vglue 0mm\line{\small\pfoot\ifnum\count0=\startpage
\copyright\ \gtp\hfill\else
\gt, Volume \thevolumenumber\ (\thevolumeyear)\hfill\fi}\vss
}}
\else
\makeatletter
\def\@oddhead{{\small\lhead\ifnum\count0=\startpage ISSN 1364-0380 (on line)
1465-3060 (printed) \hfill {\lnum\number\count0}\else\ifodd\count0
\def\\{ }\ifx\theshorttitle\relax \thetitle \else\theshorttitle\fi\hfill
{\lnum\number\count0}\else\def\\{ and }{\lnum\number\count0}
\hfill\ifx\theshortauthors\relax 
\theauthors\else\theshortauthors\fi\fi\fi}}\def\@evenhead{\@oddhead}
\def\@oddfoot{\small\lfoot\ifnum\count0=\startpage\copyright\ \gtp\hfill\else
\gt, Volume \thevolumenumber\ (\thevolumeyear)\hfill\fi}
\def\@evenfoot{\@oddfoot}
\makeatother
\fi

\newwrite\gtoutfile
\long\gdef\makeheadfile{  %%% start of definition of \makeheadfile
{\def\\{, }\def\s{ }
\immediate\openout\gtoutfile head.xxx
\immediate\write\gtoutfile{To: math@arxiv.org}
\immediate\write\gtoutfile{Subject: put or rep NNNNN:pppp}
\immediate\write\gtoutfile{--text follows this line--}
\immediate\write\gtoutfile{Proxy-for: \ifx\theasciiauthors\relax
\theauthors\else\theasciiauthors\fi\s<\ifx\theasciiemail\relax\theemail\else\theasciiemail\fi>}
\immediate\write\gtoutfile{\noexpand\\}
\immediate\write\gtoutfile{Authors: \ifx\theasciiauthors\relax
\theauthors\else\theasciiauthors\fi}
\immediate\write\gtoutfile{Title: \ifx\theasciititle\relax
\thetitle\else\theasciititle\fi}
\immediate\write\gtoutfile{Subj-class: GT or SG or MG etc}
\immediate\write\gtoutfile{MSC-class: \theprimaryclass\ifx\thesecondaryclass\relax\else, \thesecondaryclass\fi}
\immediate\write\gtoutfile{Journal-ref: Geom. Topol. \thevolumenumber
(\thevolumeyear) \startpage-\finishpage}
\immediate\write\gtoutfile{Comments: Published in Geometry and Topology at}
\immediate\write\gtoutfile{\s\s http://www.maths.warwick.ac.uk/gt/GTVol\thevolumenumber/paper\thepapernumber.abs.html}
\immediate\write\gtoutfile{\noexpand\\}
\immediate\write\gtoutfile{}
\ifx\theasciiabstract\relax
\immediate\write\gtoutfile{\theabstract}\else
\immediate\write\gtoutfile{\theasciiabstract}\fi
\immediate\write\gtoutfile{}
\immediate\write\gtoutfile{\noexpand\\}
\immediate\write\gtoutfile{}
\immediate\write\gtoutfile{<uuencoded .tar.gz file here>}
\immediate\write\gtoutfile{}
\immediate\closeout\gtoutfile}}  %%% end of definition of \makeheadfile

\def\maketitlepage{\maketitlep\makeheadfile}
\let\maketitle\maketitlepage

\volumenumber{4}\papernumber{12}\volumeyear{2000}
\pagenumbers{369}{395}
\proposed{Robion Kirby}
\seconded{Walter Neumann, David Gabai}
\received{13 April 2000}
\revised{2 November 2000}
\accepted{10 October 2000}
\published{4 November 2000}

\reflist

\key{1} {\bf D Calegari},
{\it Foliations transverse to triangulations
of 3--manifolds}, 
Comm. Analysis Geom. 8 (2000) 133--158

\key{2} {\bf D Gabai},
{\it Foliations and the topology of 3--manifolds},
J. Differ. Geom. 18 (1983) 445--503

\key{3} {\bf D Gabai},
{\it Foliations and the topology of 3--manifolds, II},
J. Differ. Geom. 26 (1987) 461--478

\key{4} {\bf M Lackenby}, 
{\it Word hyperbolic Dehn surgery}, 
Invent. Math. 140 (2000) 243--282

\key{5} {\bf U Oertel}, 
{\it Homology branched surfaces: Thurston's norm on $H_2(M^3)$}, from:
``Low-Dimensional Topology and Kleinian Groups'', (D\,B\,A Epstein
editor), London Math. Soc. Lecture Notes, 112 (1986) 253--272

\key{6} {\bf L Person}, 
{\it A piecewise linear proof that the
singular norm is the Thurston norm}, 
Top. Appl. 51 (1993) 269--289

\key{7} {\bf C Petronio}, {\bf J Porti}, 
{\it Negatively oriented ideal triangulations and a proof of Thurston's 
hyperbolic Dehn filling theorem}, preprint (1999)

\key{8} {\bf M Scharlemann}, 
{\it Sutured manifolds and generalized Thurston norms}, 
J. Differ. Geom. 29 (1989) 557--614

\key{9} {\bf A Thompson}, 
{\it Thin position and the recognition problem for $S^3$}, 
Math. Res. Lett. 1 (1994) 613--630

\key{10} {\bf W Thurston}, 
{\it The Geometry and Topology of 3--manifolds},
Princeton University (1978--79)

\endreflist

\title{Taut ideal triangulations of 3--manifolds}
\asciititle{Taut ideal triangulations of 3-manifolds}
\shorttitle{Taut ideal triangulations of 3-manifolds}
\authors{Marc Lackenby}
\address{Mathematical Institute, Oxford University
\\24--29 St Giles', Oxford OX1 3LB, UK}
\asciiaddress{Mathematical Institute, Oxford University, 
24-29 St Giles', Oxford OX1 3LB, UK}
\email{lackenby@maths.ox.ac.uk}

\asciiabstract{A taut ideal triangulation of a 3-manifold 
is a topological ideal triangulation with extra combinatorial
structure: a choice of transverse orientation on each ideal
2-simplex, satisfying two simple conditions. The aim of this paper is
to demonstrate that taut ideal triangulations are very common, and
that their behaviour is very similar to that of a taut foliation.  For
example, by studying normal surfaces in taut ideal triangulations, we
give a new proof of Gabai's result that the singular genus of a knot
in the 3-sphere is equal to its genus.}

\abstract
A taut ideal triangulation of a 3--manifold is a topological 
ideal triangulation with extra combinatorial structure: 
a choice of transverse orientation on each ideal 2--simplex, 
satisfying two simple conditions. The aim of this paper is to 
demonstrate that taut ideal triangulations are very common, 
and that their behaviour is very similar to that of a taut foliation. 
For example, by studying normal surfaces in taut ideal triangulations, 
we give a new proof of Gabai's result that the singular genus of a knot
in the 3--sphere is equal to its genus. 
\endabstract
\primaryclass{57N10}
\secondaryclass{57M25}
\keywords{Taut, ideal triangulation, foliation, singular genus}

\maketitlepage

\section{Introduction}

In his famous lecture notes [\ref{10}], Thurston introduced a surprising
topological construction of the figure-eight knot complement,
by gluing two ideal tetrahedra along their faces. Using this,
he gave the knot complement a complete hyperbolic structure.
Ideal triangulations are not only a useful tool in
hyperbolic geometry (for example [\ref{7}]), but also
provide an elegant way of visualising 3--manifolds with boundary.
In this paper, we introduce `taut' ideal triangulations,
which are ideal triangulations with a little extra structure. 
Instead of relating to hyperbolic geometry, they
are more closely associated with taut foliations. 
In his seminal paper [\ref{2}], Gabai constructed taut foliations
on many Haken 3--manifolds, via his theory of sutured manifolds.
Also utilising sutured manifolds, we will prove that many
torally bounded 3--manifolds admit a taut ideal triangulation.
An analysis of normal surfaces (and their generalisations)
in taut ideal triangulations will yield a new proof of
Gabai's result that the singular genus of a knot in $S^3$
is equal to its genus. This avoids many
of the foliation technicalities of Gabai's original
argument. We hope that taut ideal triangulations will
be useful in other areas of 3--manifold theory in
the future. Some speculations on possible applications
are included in the final section of the paper.

\rk{Definition} An {\sl ideal 3--simplex} is a 3--simplex
with its four vertices removed. An {\sl ideal triangulation}
of a 3--manifold $M$ is an expression of $M - \partial M$
as a collection of ideal 3--simplices with their faces glued in
pairs. A {\sl taut
ideal triangulation} is an ideal triangulation
with a transverse orientation assigned to each ideal 2--simplex, such
that
\items\itemb for each ideal 3--simplex, precisely two of its
faces are oriented into the 3--simplex, and precisely two
are oriented outwards, and
\itemb the faces around each edge are oriented
as shown in Figure 1: all but precisely two pairs of
adjacent faces encircling the edge have compatible orientations
around that edge.
\enditems

\figure %%1
{\epsfbox{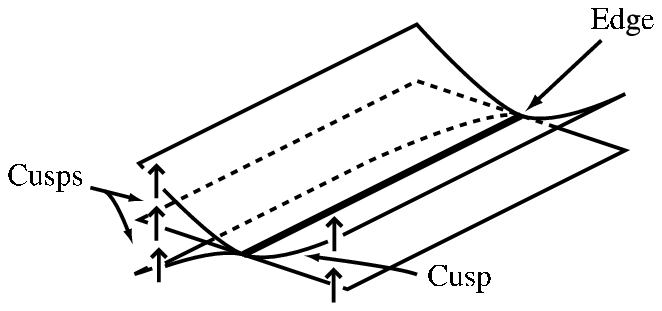}}
\endfigure

In order to describe situations such as Figure 1 more concisely,
we introduce some terminology. Suppose that some transversely
oriented surfaces embedded in a 3--manifold meet at
a 1--manifold $C$ in each of their boundaries. Then the intersection
between adjacent surfaces $S_1$ and $S_2$ is {\sl cusped} at $C$
if $S_1$ and $S_2$ are compatibly oriented around $C$.
Thus, the second of the above conditions in the definition of
a taut ideal triangulation can be rephrased as follows: all but
precisely two pairs of adjacent faces encircling an edge
have cusped intersection.

As an example, note that Thurston's ideal triangulation of the
figure-eight knot complement can be assigned a transverse orientation,
as shown in Figure 2, which makes it taut.

\figure %%2
{\epsfbox{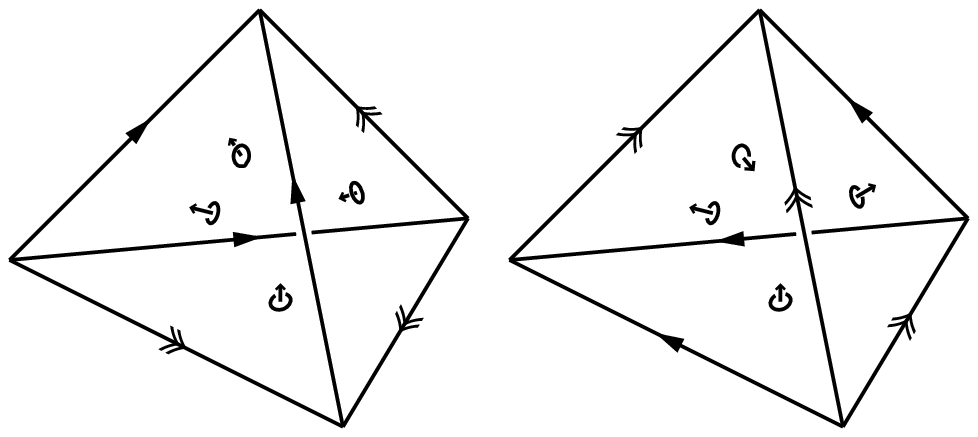}}
\endfigure

The faces and edges of a taut ideal triangulation form a
transversely oriented branched surface in $M - \partial M$.
Its branch locus is not `generic', since more than
three faces may meet at any edge. We will see that surfaces 
carried by this branched surface have strong genus-minimising 
properties. Given that taut ideal triangulations are a
special sort of branched surface, it is not surprising
that they should be related to taut foliations.
In fact, it is the absence of `triple points' in this
branched surface that gives taut ideal
triangulations many of their special properties.

The aim of this paper is to show that taut ideal
triangulations are very common, and that their presence
in a 3--manifold has useful consequences. The following is
our existence result.

\proclaim{Theorem 1} Let $M$ be a compact
orientable irreducible an-annular 3--manifold with $\partial M$
a non-empty collection of incompressible tori. Then $M$ 
has a taut ideal triangulation.\endproc

Many of the conditions in this theorem are necessary:
we will see (Proposition 10) that if a compact orientable
3--manifold admits a taut ideal triangulation, then it is
irreducible and its boundary is a non-empty collection
of incompressible tori. However, the condition that $M$ be
an-annular can be weakened.

These taut ideal triangulations are constructed from
properly embedded surfaces that are taut, in the sense
that they are incompressible and have the smallest possible Thurston complexity in
their homology class in $H_2(M, \partial M)$. Recall that
the {\sl Thurston complexity} $\chi_-(S)$ of a compact connected surface
$S$ is $\max \lbrace 0, -\chi(S) \rbrace$. The Thurston complexity
of a compact disconnected surface is defined to be the sum of
the complexities of its components. We have the following refinement of
Theorem 1, which relates taut ideal
triangulations to surfaces that minimise Thurston complexity
in their homology class.

\proclaim{Theorem 2}Let $M$ be a compact 
orientable irreducible 3--manifold with $\partial M$ a non-empty collection
of tori. Let $S$ be a properly embedded compact oriented
surface in $M$, such that
\items\itemb every component of $S$ has negative Euler characteristic
and has non-empty boundary,
\itemb $S$ has minimal Thurston complexity among all
embedded surfaces in its class in $H_2(M, \partial M)$,
\itemb for any component $T$ of $\partial M$,
the curves $T \cap \partial S$ are all essential in $T$
and coherently oriented, and
\itemb there is no properly embedded essential annulus in $M$
disjoint from $S$.
\enditems

Then $S$ is carried by the
underlying branched surface of some taut ideal
triangulation of $M$.
In fact, $S - \partial S$ is a union of ideal 2--simplices in that
ideal triangulation.\endproc

For example, the genus one Seifert surface for
the figure-eight knot is carried by the underlying
branched surface of the ideal triangulation in Figure 2.
(It lies in a regular neighbourhood of the front two faces
in the ideal 3--simplex on the right of the figure.)

We also present a converse to Theorem 2, which
asserts that surfaces carried by taut ideal triangulations
minimise Thurston complexity, even when non-embedded surfaces
are also considered.

\proclaim{\bf Theorem 3}Let $S$ be a compact
properly embedded surface carried by the underlying branched surface of
some taut ideal triangulation of a 3--manifold $M$. Then
$S$ has smallest Thurston complexity among all (possibly
non-embedded) surfaces in its class in $H_2(M, \partial M)$.\endproc

By combining Theorems 2 and 3, we retrieve Gabai's
result on the singular genus of knots. Recall that the 
{\sl singular genus} of a knot $K$ in $S^3$
is the smallest possible genus of a compact orientable surface $F$ mapped into
$S^3$ via a map $f \colon F \rightarrow S^3$ with $f^{-1}(K) = \partial F$
and $f|_{\partial F}$ an embedding onto $K$. Apply
Theorem 2, with $M$ being the exterior of $K$, and
$S$ being a minimal genus Seifert surface. Then
apply Theorem 3 to obtain the following.

\proclaim{Corollary 4}\rm [\ref{2}]\qua {\sl The singular genus
of a knot in $S^3$ is equal to its genus.}\endproc

Interestingly, Theorem 3 is proved
using normal surfaces in taut ideal triangulations.
Given that normal surfaces have a useful r\^ole to play
in other areas of 3--manifold topology (for example [\ref{9}]),
this suggests that taut ideal triangulations will have
other interesting applications. 

Taut ideal triangulations are closely related to angled ideal
triangulations, defined and studied by Casson, and developed in [\ref{4}]. 
An {\sl angled ideal triangulation} is an ideal triangulation with
a number in the range $(0, \pi)$ assigned to each
edge of each ideal 3--simplex, known as the interior
angle at that edge. These angles are required to satisfy
two simple conditions: the angles around an edge sum to
$2\pi$; and the angles at each ideal vertex of each
ideal 3--simplex sum to $\pi$. Taut
ideal triangulations induce a similar structure, except
that there are only two options for the interior angles:
the cusped intersections between faces have zero interior angle
and the non-cusped intersections have interior
angle $\pi$. (See Figure 3.) 

\figure %%3
{\epsfbox{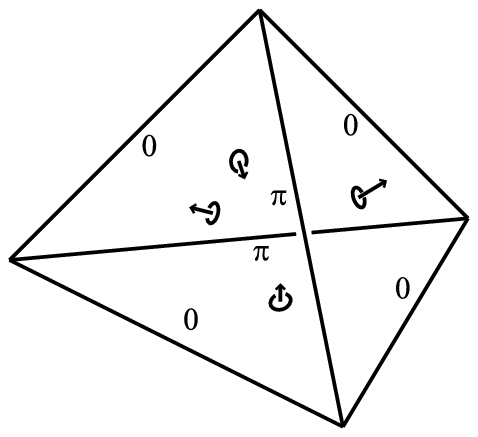}}
\endfigure

The underlying branched surface
of a taut ideal triangulation carries an essential lamination.
This is formed by laminating a neighbourhood of 
each face using a Cantor set transversal, and then
patching these laminations together at the edges.
Since the branched surface has no vertices, there
is no obstruction to performing this procedure.
This lamination extends to a foliation of the 3--manifold $M$, since the
complimentary regions of the branched surface are
products. It is not hard to find a closed curve transverse
to the foliation intersecting every leaf, and so the foliation
is taut. Note also that if $S$ is any surface carried
by the branched surface, then we may decompose
a neighbourhood $N$ of the branched surface along $S$,
then laminate $N - S$, and then extend this to a taut foliation
of $M$ in which $S$ is a leaf. Hence, by a theorem of
Gabai [\ref{2}], $S$ must have smallest Thurston complexity
among all (possibly non-embedded) surfaces in its class
in $H_2(M, \partial M)$. We therefore obtain Theorem 3.
However, one of the aims of this paper is, of course,
to provide a proof of Theorem 3 using an argument
that avoids foliations.

The underlying branched surface of a taut ideal triangulation
is, in the terminology of [\ref{5}], an example of a taut homology
Reebless incompressible branched surface. In [\ref{5}], Oertel
analysed this type of branched surface and used them to
establish certain facts about the Thurston norm. However,
they do not seem to interact so well
with singular surfaces as do our taut ideal triangulations.

A purely combinatorial proof of Corollary 4 (and, more generally,
of the equivalence of the Thurston norm and the singular norm)
has been given by Person in [\ref{6}], building
on Scharlemann's foliation-free approach to sutured manifolds [\ref{8}].
The argument in [\ref{6}] is rather different from the proof
given here. In particular, it does not deal with ideal triangulations or
normal surfaces.

\section{Constructing taut ideal triangulations}

In this section, we prove Theorems 1 and 2, by constructing
the required taut ideal triangulations. Suppose therefore
that $M$ is a compact orientable irreducible an-annular
3--manifold, with $\partial M$ a non-empty union of  incompressible tori.
The simplest case is where $M$ fibres over the circle with fibre $S$, and
then the result is rather easy. Note that $S$ has negative Euler characteristic
and non-empty boundary, by our assumptions on $M$, and
hence $S$ has an ideal triangulation. Using this, we will construct
a taut ideal triangulation of $M$.
We therefore now recall some well-known facts about
ideal triangulations of surfaces.

It will be helpful to consider the following generalisation
of an ideal triangulation.
An {\sl ideal region} $\delta$ for a compact
surface $S$ is a compact submanifold of
$\partial S$ having non-empty intersection with
each component of $\partial S$. An {\sl ideal triangulation}
of $S$ with ideal region $\delta$ is an expression of
$S - \delta$ as a union of ideal 2--simplices with some of
their edges glued in pairs. Hence, $\partial S -\delta$ must be a
(possibly empty) collection of open
arcs, each of which is an edge of an ideal 2--simplex.
We define the {\sl triangular number} $t(S, \delta)$ of a surface $S$
with ideal region $\delta$ to be
$$t(S, \delta) = -2\chi(S) + |\partial S - \delta|.$$

\proclaim{Lemma 5}Let $S$ be a compact orientable surface
with non-empty boundary and ideal region $\delta$.
If $t(S, \delta) > 0$, then $S$ admits an ideal triangulation
with ideal region $\delta$. Any such ideal triangulation
contains precisely $t(S, \delta)$ ideal 2--simplices.\endproc

\prf The surface $S - \delta$ is obtained from
a compact orientable surface $\hat S$ by removing a finite number of points $P$
from its boundary and from its interior. If $\hat S$
is closed and is not a 2--sphere, it has a one-vertex
triangulation. If $\hat S$ has non-empty boundary and
is not a disc, it has a triangulation with a single vertex
on each boundary component and no vertices in its interior.
Subdivide these triangulations, if necessary,
so that its vertices are precisely $P$.
Then remove these vertices to obtain an ideal triangulation
of $S - \delta$. The argument when $\hat S$ is a sphere or
a disc is similar. We find a triangulation of $\hat S$, and then
remove its vertices to obtain an ideal triangulation of
$S - \delta$. The assumption that $t(S, \delta) > 0$ guarantees that
this is possible.

Now consider any such ideal triangulation of $S$. It is
formed by gluing the edges of ideal triangles in pairs. Each ideal
triangle has triangular number one. At each gluing of edges, the total
Euler characteristic goes down by one, but the number
of boundary edges goes down by two. Hence, the total
triangular number is unchanged. Thus, this ideal
triangulation has $t(S, \delta)$ ideal 2--simplices.\break
\hbox to 5pt{\hss} \endprf

It is also very well known that any two ideal triangulations
of $S$ with the same ideal region
differ by a sequence of the following {\sl elementary
moves}: pick two distinct ideal 2--simplices that share
an edge; remove this edge, forming a `square' whose interior
is embedded in the interior of $S$; then subdivide this
square along its other diagonal to form two new ideal
triangles. This fact is so important for our presentation
that we include a proof.

\figure %%4
{\epsfbox{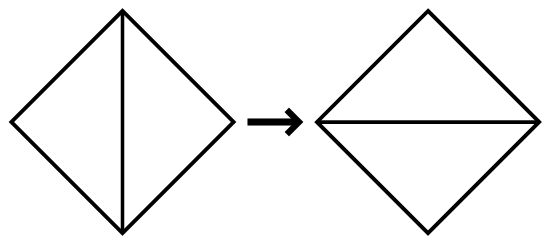}}
\endfigure

\proclaim{Lemma 6}Any two ideal triangulations
of a compact orientable surface $S$ with the same ideal region
$\delta$ are related by a finite sequence of
elementary moves and an ambient isotopy.\endproc

\prf Let $T_1$ and $T_2$ be the ideal
triangulations of $S$. Note that $T_1$ and $T_2$
both have $t(S, \delta)$ ideal triangles. Let $E_2$ be the edges of $T_2$.
We will work with the dual graph $G_1$ of $T_1$, which
is a graph embedded in $S$, the interior vertices having
valence three, and vertices on each component of $\partial S - \delta$
having valence one. 

Perform a small ambient isotopy so that the intersection
between $G_1$ and $E_2$ is transverse and disjoint from
the trivalent vertices of $G_1$. Note that
each component of $(S  - \delta) - G_1$ is either 
$I \times {\Bbb R}$ or $S^1 \times {\Bbb R}$. We consider three possibilities.

{\bf Case 1}\qua Some edge of $E_2$ is disjoint
from $G_1$.

In this case, the arc lies entirely in a component 
$I \times {\Bbb R}$ or $S^1 \times {\Bbb R}$ of
$(S - \delta) - G_1$, and both its ends lie in the
same end of $I \times {\Bbb R}$ or $S^1 \times {\Bbb R}$.
An extrememost such arc separates
off a disc of $S - E_2$ with a single end. 
However, every component of $S - E_2$
is a triangle, and we therefore have a contradiction.
Thus, this case does not arise.

{\bf Case 2}\qua Every arc of $E_2 - G_1$
runs between $G_1$ and an end of $S - \delta$.

Then each edge of $E_2$ intersects $G_1$ at a single point.
Hence, each triangle of $T_2$ has three arcs of $G_1$ entering
it. It therefore has at least one trivalent vertex of $G_1$
in its interior. However,
there are as many trivalent vertices of $G_1$ as there are triangles
of $T_2$, and so each triangle of $T_2$ contains a single trivalent
vertex of $G_1$. Hence, $G_1$ is the dual of $T_2$,
and so (up to ambient isotopy) $T_1$ and $T_2$ are
the same ideal triangulation.

{\bf Case 3}\qua Some arc of $E_2 - G_1$ has
both endpoints in $G_1$.

Pick such an arc $\alpha$ extrememost in one component
of $S - G_1$. This separates off a disc $D$ which
lies in some triangle of $T_2$. Let $\beta$
be $\partial D - \alpha$, which is a path in $G_1$.
If $\beta$ runs through at most one vertex of $G_1$, then there
is an ambient isotopy of $G_1$ which reduces
$\vert G_1 \cap E_2 \vert$. Suppose therefore
that $\beta$ contains at least two vertices of $G_1$.
Apply an elementary move to adjacent vertices of
$G_1 \cap \beta$ to reduce the number of vertices of
$\beta$. Repeat this until $\beta$ contains only one vertex,
and then perform an ambient isotopy to reduce
$\vert G_1 \cap E_2 \vert$. In this
way, we remove all arcs of $E_2 - G_1$ with
both endpoints in $G_1$. Hence, after a finite
number of elementary moves, we end with $T_1$
ambient isotopic to $T_2$. \endprf

Let us now return to the manifold $M$ which fibres over
the circle, with fibre $S$. 
Let $h \colon S \rightarrow S$ be the monodromy homeomorphism. 
Pick some ideal triangulation $T$ of $S$ with ideal region $\partial S$.
Then $h(T)$ may be transformed into 
$T$ by sequence of elementary moves and an
isotopy. Of course, in this sequence of moves, we may guarantee
that no edge of $T$ is left untouched. (For example, if an edge
is adjacent to two distinct 2--simplices, 
perform an elementary move and its reverse).
The taut ideal triangulation of $M$ is constructed
as follows. Start with $S$ and its triangulation $T$.
Each time that an elementary move is performed, glue an
ideal tetrahedron (as in Figure 5) onto one side of $S$,
attaching it to the
two ideal triangles involved in the elementary move.
This side of $S$ then inherits the new ideal triangulation.
Repeat this process for each elementary move. Then we have
constructed $(S - \partial S) \times I$, since every
edge of $T$ was modified by the elementary moves.
Now glue the two components of $(S - \partial S) \times \partial I$
via $h$. The ideal triangulations match up to form a taut
ideal triangulation of $M$.
The choice of elementary moves realizing $h$ was highly non-unique, and
hence $M$ has many taut ideal triangulations.
It is interesting to note that the taut
ideal triangulation of the figure-eight knot
complement given in Figure 2 can be constructed in this way,
except that only two elementary moves are used, and hence
one edge of the ideal triangulation of the fibre remains unmodified.

\figure %%5
{\epsfbox{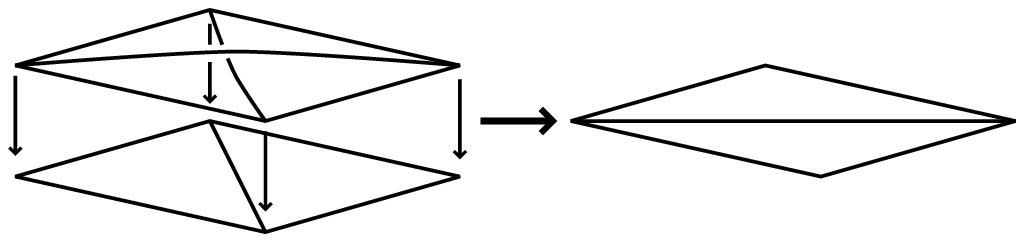}}
\endfigure

When $M$ does not fibre over $S^1$, it is significantly
more difficult to construct taut ideal triangulations. 
The main technical tool is sutured manifold theory, which Gabai
originally used to construct taut foliations on many irreducible
3--manifolds. We will use the version developed by
Scharlemann in [\ref{8}], and, in this section of the paper,
we will assume that the reader is reasonably familiar with [\ref{8}].
Recall that a {\sl sutured manifold} $(M, \gamma)$
is a compact oriented 3--manifold $M$, with $\partial M$
decomposed into subsurfaces ${\cal R}_-$, ${\cal R}_+$,
$A(\gamma)$ and $T(\gamma)$, which intersect in simple closed curves, such that
\items \itemb each component of $A(\gamma)$ is an annulus
adjacent to both ${\cal R}_-$ and ${\cal R}_+$,
\itemb each component of $T(\gamma)$ is a torus, and
\itemb ${\cal R}_- \cap {\cal R}_+ = \emptyset$.
\enditems

We let $\gamma$ be the core curves of $A(\gamma)$.
The surface ${\cal R}_-$ (respectively, ${\cal R}_+$)
is assigned an orientation into (respectively, out of)
$M$. The annuli $A(\gamma)$ and tori $T(\gamma)$ are not assigned a specific
orientation.

When $S$ is a transversely oriented surface properly embedded
in a sutured manifold $(M, \gamma)$, the transverse orientations on $S$
and ${\cal R}_\pm$ induce a cusp on one side of each component of $\partial S
\cap {\cal R}_\pm$. The manifold $M - {\rm int}({\cal N}(S))$
inherits a sutured manifold structure, providing $S$ satisfies various simple
properties. These properties have a variety of names in the
literature: d--surface or conditioned in [\ref{8}], groomed or well-groomed in [\ref{3}].
In this paper, we introduce a variant of these. We
allow $S$ to intersect $A(\gamma)$ in simple closed curves,
{\sl transverse} arcs (which run between distinct components
of $\partial A(\gamma)$) and {\sl glancing} arcs (which
run between the same component of $\partial A(\gamma)$).

\rk{Definition} A transversely oriented surface $S$ properly
embedded in $(M, \gamma)$ is {\sl styled} if
\items \itemb for each component $T$ of $T(\gamma)$,
the curves $T \cap \partial S$ are all essential in $T$ and coherently oriented,
\itemb near each simple closed curve of $S \cap A(\gamma)$,
$S$ has the same orientation as ${\cal R}_-$ and ${\cal R}_+$ near
that component of $A(\gamma)$,
\itemb the transverse arcs of intersection between $S$ and
any component of $A(\gamma)$ are all coherently oriented,
and
\itemb near any glancing arc $\alpha$ of $S \cap A(\gamma)$,
the cusped side of $\partial S$ runs into the disc component of
$A(\gamma) - \alpha$.
\enditems

\figure %%6
{\epsfbox{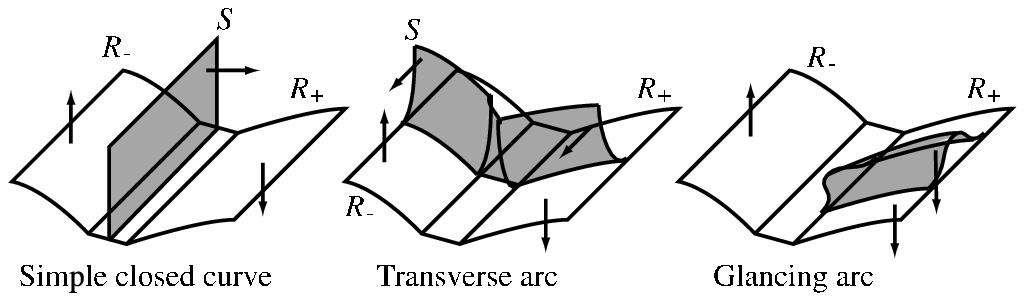}}
\endfigure

(See Figure 6.) When $S$ is styled, then $M_S = M - {\rm int}({\cal N}(S))$ inherits
a sutured manifold structure $(M_S,\gamma_S)$ in a natural way:
${\cal R}_\pm(\gamma_S)$ is composed of ${\cal R}_\pm(\gamma) - {\rm int}({\cal N}(S))$
and the copies of $S$. The tori of $T(\gamma)$ disjoint from $S$
yield $T(\gamma_S)$. The annuli $A(\gamma_S)$ lie between
${\cal R}_-(\gamma_S)$ and ${\cal R}_+(\gamma_S)$, and
also arise from $T(\gamma) - {\rm int}({\cal N}(S))$. The
orientation conditions imposed on glancing arcs and simple
closed curves of $S \cap A(\gamma)$
guarantee that each component of $A(\gamma) - {\rm int}({\cal N}(S))$ lies in
$A(\gamma_S)$.

We introduce the following definition.

\rk{Definition} An {\sl ideal region} $\delta$ of a sutured manifold
$(M, \gamma)$ is a collection of the following subsets of
$A(\gamma) \cup T(\gamma)$:
\items \itemb all of $T(\gamma)$,
\itemb possibly some components of $A(\gamma)$, and
\itemb `squares', each of which is a subset
of a component $A$ of $A(\gamma)$, being the region between
two properly embedded transverse arcs in $A$.
\enditems

\nobreak
We insist, in addition, that no component of $A(\gamma)$ is 
disjoint from $\delta$.

The idea behind the above definition is that we start with 
a sutured manifold $(M, \gamma)$ having $\partial M = T(\gamma)$,
and so the whole of $\partial M$ is the ideal region.
Then we perform a sequence of sutured manifold decompositions,
resulting in sutured manifolds embedded in $M$. Their
ideal regions will be their intersection with $\partial M$.

The taut ideal triangulations in Theorems 1
and 2 will be constructed using a sutured manifold hierarchy.
At each stage of this hierarchy, we will construct a
taut triangulation, in the following sense.

\rk{Definition} A {\sl taut triangulation}
of a sutured manifold $(M, \gamma)$ with ideal region
$\delta$ is an expression of $M - \delta$ as
a collection of ideal 3--simplices with some of their
faces identified in pairs, and with a transverse orientation 
assigned to each ideal 2--simplex, such that
\items \itemb for each ideal 3--simplex, precisely two of its
faces are oriented into the 3--simplex, and precisely two
are oriented outwards,
\itemb each component of $\gamma - \delta$ is an
edge of the triangulation, 
\itemb each 2--simplex in $\partial M$ lies entirely
in ${\cal R}_-$ or ${\cal R}_+$ (apart from a collar neighbourhood
of some of its edges, which may lie in $A(\gamma)$),
\itemb the
transverse orientation of each 2--simplex in $\partial M$
agrees with that of ${\cal R}_\pm$,
\items \itemb for each edge not in $\gamma$, 
all but precisely two pairs of adjacent faces around that edge 
have cusped intersection, and
\itemb for each edge in $\gamma$, all faces
around that edge have cusped intersection.
\enditems

\figure %%7
{\epsfbox{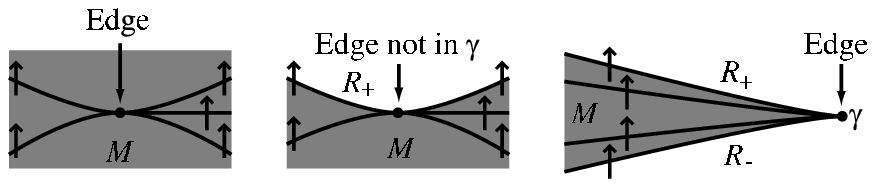}}
\endfigure

When $(M, \gamma)$ has a taut
triangulation, ${\cal R}_\pm$ inherits an ideal triangulation
with ideal region $\delta \cap {\cal R}_\pm$.
Note also that when $\partial M$ is a collection of tori, and
$\delta = T(\gamma)$ is all of these tori, then a taut
triangulation of $M$ is a taut ideal triangulation.

The case where $M$ is a 3--ball and $\gamma$ is a single
curve is an instructive example.
If $\delta$ is precisely four squares, then a single
ideal 3--simplex forms a taut triangulation of $M$, as in
Figure 3. If $\delta$ is more than four squares, 
then we may pick an ideal triangulation
of ${\cal R}_- - \delta$.
Then, as in the fibred case, a taut 
triangulation of $M$ arises by attaching a collection
of ideal 3--simplices to realize a suitable sequence of 
elementary moves. However, if $\delta$ is fewer
than four squares or the whole of $A(\gamma)$, then 
$t({\cal R}_-, \delta)$ and $t({\cal R}_+, \delta)$ are
each at most one. Thus, any ideal triangulation
of ${\cal R}_- - \delta$ and ${\cal R}_+ -
\delta$ consists of at most one ideal triangle.
It is then not hard to see that $(M,\gamma)$ admits
no taut triangulation having $\delta$ as an ideal region.
The following theorem demonstrates that 
a similar pattern arises for other taut sutured manifolds.
We say that an annulus embedded in $M$ is {\sl $\delta$--essential}
if
\items \itemb its boundary is in $\delta$,
\itemb it is incompressible, and
\itemb it is not parallel to an annulus
in $\delta$.
\enditems

\proclaim{Theorem 7}Let $(M, \gamma)$ be a
connected taut sutured manifold. Suppose that 
$\partial M$ is non-empty and that no component
of $\partial M$ is disjoint from $A(\gamma) \cup T(\gamma)$.
Let $\delta$ be any choice of ideal region,
such that $t({\cal R}_-, \delta) \geq 2$ and $t({\cal R}_+, \delta) \geq 2$,
and where $M$ contains no $\delta$--essential annulus.
Then $(M, \gamma)$ admits a taut triangulation with $\delta$ as ideal
region.\endproc

{\bf Proof of Theorems 1 and 2 from Theorem 7}\qua
We first show that Theorem 1 follows from Theorem 2.
Let $M$ be a compact orientable irreducible an-annular 3--manifold with
$\partial M$ a non-empty collection of incompressible tori. 
Then, it is well known that some element of
$H_2(M, \partial M)$ has non-trivial image in 
$H_1(\partial M)$ under the boundary map of the
homology exact sequence of the pair $(M, \partial M)$.
This element of $H_2(M, \partial M)$ is
represented by a properly embedded compact
orientable surface $S$ with non-empty boundary.
Take $S$ to have minimal Thurston complexity among
all embedded surfaces in its homology class.
We may cap off curves of $\partial S$ bounding
discs in $\partial M$, and also, by attaching annuli,
we may assume that, for each component $T$ of $\partial M$,
$\partial S \cap T$ is a collection of coherently
oriented simple closed curves, each essential in $T$.
This does not increase its Thurston complexity or
alter its class in $H_2(M, \partial M)$. 
By construction, at least one component of $S$ has non-empty
boundary. Restrict attention to this component, which we
will now call $S$. Then, $S$ is neither a disc nor an annulus,
since $M$ has incompressible boundary and is an-annular.
Hence, $S$ satisfies all of the conditions
of Theorem 2, and so Theorem 1 follows
from Theorem 2. 

We now show that Theorem 2 follows from Theorem 7.
Let $M$ and $S$ be as in Theorem 2, and
let $T(\gamma) = \partial M$. Then $(M, \gamma)$ is taut
since $M$ is irreducible and not a solid torus.
Perform the taut sutured manifold decomposition
$$(M, \gamma) \buildrel S \over \longrightarrow (M_S, \gamma_S),$$
and let $\delta_S = A(\gamma_S) \cup T(\gamma_S)$.
Note that there is no $\delta_S$--essential annulus
in $M_S$, since this would be an essential annulus
in $M$ disjoint from $S$. Also,
$$t({\cal R}_-(\gamma_S), \delta_S) =
t({\cal R}_+(\gamma_S), \delta_S) = - 2 \chi(S) \geq 2.$$
Thus, using Theorem 7, find a taut triangulation
for $(M_S, \gamma_S)$ with ideal region $\delta_S$.
Glue the two copies of $S$ in ${\cal R}_-(\gamma_S)$
and ${\cal R}_+(\gamma_S)$ together, ensuring that their ideal triangulations
agree using a sequence of elementary moves, as in the
fibred case. The result is a taut ideal triangulation
of $M$ in which $S$ is a union of ideal 2--simplices
and is carried by the underlying branched surface. \endprf

The remainder of this section is devoted to proving
Theorem 7. Let $(M, \gamma)$ and $\delta$ be as in Theorem 7.
We will prove the theorem by induction
backwards along a sutured manifold hierarchy for $(M, \gamma)$.
Each decomposing surface (other than product discs)
will end up as a collection of ideal 2--simplices
in the resulting taut triangulation. Hence
it is vital that each surface meets $\delta$. We do
this by `sliding' the boundary of the surface towards 
$\delta$ along arcs. The following lemma guarantees that
this is possible.

\proclaim{Lemma 8}Let $(M,\gamma)$
be a connected taut sutured manifold where $\partial M$ is non-empty
and no component of
$\partial M$ is disjoint from $A(\gamma) \cup T(\gamma)$. 
Then there is a taut sutured manifold hierarchy
$$(M, \gamma) = (M_1, \gamma_1) 
\buildrel S_1 \over \longrightarrow
(M_2, \gamma_2) 
\buildrel S_2 \over \longrightarrow \dots
\buildrel S_{n-1} \over \longrightarrow
(M_{n}, \gamma_{n}),$$
such that, for each $i$, $S_i$ is a connected styled non-separating
surface with non-empty boundary, and,
for any point $p$ on ${\cal R}_\pm(\gamma_i)$,
there is an embedded arc $\alpha$ in $\partial M_i$
such that
\items \itemb $p$ is one endpoint of $\alpha$,
\itemb $\alpha \cap A(\gamma_i)$ is the other endpoint of $\alpha$,
and
\itemb if $\alpha$ is oriented from $p$
to $\partial \alpha - p$, then at each point of intersection between $\alpha$
and $\partial S_i$, $\alpha$ runs from the cusped side of $\partial S_i$
to the uncusped side.
\enditems\endproc

\prf We must return to the proof of 
the existence of sutured manifold hierarchies (Theorem 4.19
of [\ref{8}]). There are two main ingredients: to show
that given any non-zero class $z \in H_2(M_i, \partial M_i)$,
we can perform a taut decomposition along a styled surface
$S_i$ with $z = [S_i, \partial S_i]$; then to show
that, with the correct choice of decomposing surfaces,
a sequence of taut decompositions can be made to terminate in a collection of
3--balls. The second part is dealt with in \S4 of [\ref{8}]. There,
a complexity for a sutured manifold is defined.
If a taut sutured manifold contains a non-trivial product
disc, then decomposing along this disc does not increase 
the complexity. (A product disc is {\sl non-trivial}
if it does not separate off a 3--ball.)
If a taut sutured manifold contains no non-trivial product
disc, then any taut decomposition along a connected non-separating 
incompressible surface decreases sutured manifold complexity
(Theorem 4.17 of [\ref{8}]). If a sutured manifold is decomposed
along a product disc, then future decomposing surfaces
can be ambient isotoped so that they avoid the
two copies of this disc. Hence, by Lemma 4.2 of [\ref{8}],
it is possible to postpone all the decompositions
along product discs until the final step. Hence, 
a sequence of taut decompositions along connected non-separating
incompressible surfaces must eventually terminate with a product
sutured manifold. Decompose this along non-separating 
product discs to obtain a 3--ball.
Hence, providing at each stage
we can find a surface satisfying the requirements
of Lemma 8, this sequence of sutured manifolds
can be guaranteed to terminate.

Suppose therefore that we have constructed
the sutured manifold sequence as far as $(M_i, \gamma_i)$.
We claim that no component $F$ of $\partial M_i$ is disjoint
from $A(\gamma_i) \cup T(\gamma_i)$. Suppose that, on the
contrary, $F$ is disjoint from $A(\gamma_i) \cup T(\gamma_i)$.
We are assuming that no component of $\partial M$ is
disjoint from $A(\gamma) \cup T(\gamma)$, and hence
$F$ must intersect some $S_j$, $j <i$. Let $j$ be the largest
such integer. Note that $\partial S_j$ has non-empty boundary.
If $\partial S_j$ intersects $A(\gamma_j) \cup T(\gamma_j)$,
then $F \cap A(\gamma_i)$ is non-empty. If $\partial S_j$ 
is disjoint from $A(\gamma_j) \cup T(\gamma_j)$, then
there is an arc $\alpha$ as in the lemma, which runs
from $\partial S_j$ to $A(\gamma_j)$. Whether or not the interior of
$\alpha$ intersects $\partial S_j$, we obtain
a component of $A(\gamma_i)$ in $F$. 

Suppose that $\partial M_i$ is not a collection of 2--spheres.
For otherwise we have constructed the required hierarchy.
Let $C$ be a finite collection of disjoint oriented
simple closed curves in $\partial M_i$ satisfying the following:
\items \itemb $[C] \not=0 \in H_1(\partial M_i)$,
\itemb $[C] = 0 \in H_1(M_i)$,
\itemb $C$ intersects $A(\gamma_i)$ only in transverse arcs, and
\itemb $|C|$ is minimal among all oriented
curves in $\partial M_i$ with the above three properties.
\enditems

The existence of $C$ is a consequence of the well-known
fact that $H_1(\partial M_i) \rightarrow H_1(M_i)$ has non-zero
kernel. The orientation on $C$ and some orientation on
$\partial M_i$ induce a transverse orientation on $C$.
Our aim is to construct a taut decomposing surface $S_i$
such that, away from a regular neighbourhood of $A(\gamma_i) \cup T(\gamma_i)$,
$\partial S_i$ agrees with $C$. Hence, we now check that for any point $p$ on 
${\cal R}_\pm(\gamma_i)$, we may find a path
$\alpha$ as in the lemma.

Construct a graph in $\partial M_i$, with a single vertex
in each component of $\partial M_i - C$ and with an
edge for each component of $C$ intersecting that
component transversely and missing all other components of $C$.
Orient the edges according to the transverse orientation of $C$.
Note that each vertex of the graph has valence more
than one. Otherwise, there would be a separating component of $C$, 
which would contradict the minimality assumption on
$|C|$. Note also that no vertex
of the graph can have more than one edge entering it,
or more than edge leaving it. For, in this case, these
edges correspond to distinct components of
$C$, which we may join by an arc in $\partial M_i - C$, and using this arc, we may
reduce the number of components of $C$. This contradicts the
minimality assumption again. Hence, each vertex of the 
graph has precisely one edge entering
and one edge leaving. Therefore, the graph is a disjoint
union of circles. Let $F$
be the component of $\partial M_i$ containing $p$.
The component of $F - C$ containing $p$
corresponds to a vertex $v_1$ in the graph. Some component of $F - C$
must intersect $\gamma_i$. This corresponds to a
vertex $v_2$ of the graph. There are paths in the graph from $v_1$
to $v_2$ that are compatible and incompatible with the orientation on the graph.
Truncate these paths at their first intersection points
with $A(\gamma_i)$. One of these paths is the path
$\alpha$ as required.

Transversely orient the curves $\gamma_i$ so that they point towards
${\cal R}_+(\gamma_i)$. Let $C'$ be the double-curve sum
of $C$ with a sufficient number of parallel copies of $\gamma_i$,
so that, after a small ambient isotopy, $C'$ intersects
any component of $A(\gamma_i)$ in a collection of 
coherently oriented transverse arcs.
Apply Theorem 2.5 of [\ref{8}] to these curves $C'$. This results in
a surface $S_i$ properly embedded in $M_i$, such that
$(M_i, \gamma_i) \buildrel S_i \over \longrightarrow
(M_{i+1}, \gamma_{i+1})$ is taut, and such that
$\partial S_i - A(\gamma_i)$ and $C' - A(\gamma_i)$
are the same 1--manifolds with the same transverse orientations. 
After possibly capping off oppositely oriented
simple closed curves of $\partial S_i$ in $A(\gamma_i)$
and $T(\gamma_i)$, $S_i$ becomes styled. For any point $p$ on ${\cal R}_\pm(\gamma_i)$, 
we may find an arc $\alpha$ running from $p$ to $A(\gamma_i)$
intersecting $\partial S_i$ correctly. The particular choice of transverse
orientation on $\gamma_i$ guarantees that this is
still true at the intersection points between
$\alpha$ and $\partial S_i$ coming from the parts of $\gamma_i$ in $C'$.
By restricting to some component
of $S_i$, we may ensure that $S_i$ is connected and non-separating,
and has non-empty boundary.
\endprf

{\bf Proof of Theorem 7}\qua Let
$(M, \gamma)$ and $\delta$ be as in Theorem 7. 
We consider a sutured manifold hierarchy as in Lemma 8,
and prove the theorem by induction backwards along the
hierarchy. The hierarchy ends with
a 3--ball. We have already shown in this
case that Theorem 7 holds.
We now prove the inductive step. Consider
a taut sutured manifold decomposition
$$(M, \gamma) \buildrel S \over \longrightarrow
(M_S, \gamma_S)$$
where $S$ is a surface satisfying the requirements
of Lemma 8. We assume inductively that $(M_S, \gamma_S)$
satisfies the conclusion of the theorem. 

Note that $S$ has non-negative triangular number. For, otherwise,
it is a disc intersecting $\delta$ in at most one arc (which
must be glancing) or in a simple closed curve. Since
$S$ is non-separating,
${\cal R}_\pm(\gamma)$ is compressible, which is a contradiction.

We now perform a sequence of ambient isotopies
to $S$ so that, afterwards
\items \itemb $S$ remains styled,
\itemb no component of $\partial S$ is disjoint
from $\delta$,
\itemb each arc of intersection between $S$ and $A(\gamma)$
lies in $\delta$ (but closed curves of $S \cap A(\gamma)$ need not lie
wholly in $\delta$), and
\itemb $(M_S, \gamma_S)$ remains unchanged, up to
homeomorphism.
\enditems

Note that when the above conditions hold,
the manifold $(M_S, \gamma_S)$ inherits the
ideal region $\delta_S = \delta \cap (A(\gamma_S) \cup T(\gamma_S))$. Note in particular 
that no component of $A(\gamma_S)$ is disjoint from $\delta_S$.
Also, $M_S$ contains no $\delta_S$--essential annulus,
since this would be a $\delta$--essential annulus in $M$.

We may clearly perform an ambient isotopy of $S$, supported in
a neighbourhood of $A(\gamma)$ to ensure that each arc of
intersection between $S$ and $A(\gamma)$ lies in $\delta$. 
Suppose that some component of $\partial S$
is disjoint from $A(\gamma) \cup T(\gamma)$. Pick a point $p$ on this
component. Let $\alpha$ be the arc running from $p$
to $A(\gamma)$, as in Lemma 8. We may assume that $\alpha$
avoids all components of $\partial S$ intersecting $A(\gamma)$,
and that the endpoint of $\alpha$ is in $\delta$. Then
let $p'$ be the point of $\alpha \cap S$ closest to
$A(\gamma)$. Ambient isotope the component of $\partial S$
containing $p'$ along $\alpha$ so that afterwards it does intersect $A(\gamma)$
in a glancing arc. Note that the assumption in Lemma 8
that at each intersection point with $\partial S$, $\alpha$ runs from the cusped side of
$\partial S$ to the uncusped side guarantees that $S$ is styled and that $(M_S, \gamma_S)$
remains unchanged. Repeat as necessary, until $S$ satisfies the above four conditions.

Note that $M_S$ is connected, since $S$ is connected and
non-separating. We have that
$$t({\cal R}_+(\gamma_S), \delta_S) = t({\cal R}_+(\gamma), \delta) + t(S,\delta) 
\geq t({\cal R}_+(\gamma), \delta)
\geq 2,$$ and a similar inequality holds for $t({\cal R}_-(\gamma_S), \delta_S)$. 

Our aim now is to alter $S$ so that afterwards each
component of ${\cal R}_\pm(\gamma) - {\rm int}({\cal N}(S))$
has non-negative triangular number.
A component $D$ of ${\cal R}_\pm(\gamma) - {\rm int}({\cal N}(S))$ with 
negative triangular number
must be a disc intersecting $\delta$ in at most one arc or simple closed
curve. It cannot be a component of ${\cal R}_\pm(\gamma)$, since
that would imply that $M$ was a 3--ball with
$t({\cal R}_-(\gamma), \delta) = t({\cal R}_+(\gamma), \delta) \leq -1$.
Hence $D$ must intersect $\partial S$ in a single arc.
The cusped side of $\partial S$ cannot lie in $D$, for
otherwise ${\cal R}_\pm(\gamma_S)$ would be a disc
and so $(M_S, \gamma_S)$ would be a 3--ball with
$t({\cal R}_-(\gamma_S), \delta_S) = t({\cal R}_+(\gamma_S), \delta_S) = -
1$, which is a contradiction.
Ambient isotope this arc $\partial D \cap S$ across $D$ into $A(\gamma)$ 
to reduce the number of components of $\partial S - \delta$. 
It is straightforward to check that $S$
remains styled. Hence, eventually, each
component of ${\cal R}_\pm(\gamma) - {\rm int}({\cal N}(S))$
has non-negative triangular number. 

The components of ${\cal R}_\pm(\gamma) - {\rm int}({\cal N}(S))$
with zero triangular number are of two possible types:
\items 
\item{(i)} a parallelity region between an arc $C$ of 
$\partial S \cap {\cal R}_\pm(\gamma)$ and an arc of $\gamma - \delta$, or
\item{(ii)} a region that lies between two parallel arcs $C_1$ and $C_2$
of $\partial S \cap {\cal R}_\pm(\gamma)$ and that misses $\gamma - \delta$.
\enditems

In case (i), note that the cusped side of
$C$ cannot lie in the parallelity region.
For this would create a disc component of
${\cal R}_\pm(\gamma_S)$ with zero triangular number.
Then $(M_S, \gamma_S)$ would be a 3--ball with
$t({\cal R}_-(\gamma_S), \delta_S) = t({\cal R}_+(\gamma_S), \delta_S) =
0$, which is a contradiction.
We examine the arcs of $\partial S \cap A(\gamma)$ adjacent to
$C$. They cannot both be transverse arcs, since $S$ is
styled. If they are the same glancing arc, then
we ambient isotope $C$ into $A(\gamma)$. Otherwise,
we perform an ambient isotopy of $S$ which removes a
glancing arc of $S \cap A(\gamma)$. Note that $S$ remains
styled.

In case (ii), this parallelity region cannot
contain the cusped sides of both $C_1$ and $C_2$.
For, again, this would create a disc component of
${\cal R}_\pm(\gamma_S)$ with zero triangular number.
We perform a small homotopy which amalgamates these parallel
arcs into one. Of course, this renders $S$ no longer
embedded. Note that the only possible obstruction to performing
all these homotopies is when the type (ii) parallelity regions
lie in an annular component of ${\cal R}_\pm(\gamma)$ with its entire 
boundary in $\delta$.
However, if such an annulus existed, we could push its interior
a little into the interior of $M$, forming a $\delta$--essential
annulus, which is a contradiction.

These alterations to $S$
will not alter the homeomorphism type of $(M_S, \gamma_S)$.
They may alter $\delta_S$, but both
$t({\cal R}_-(\gamma_S), \delta_S)$ and $t({\cal R}_+(\gamma_S), \delta_S)$
remain at least two. Inductively, therefore, $(M_S,\gamma_S)$ has a taut triangulation
with ideal region $\delta_S$. We now use this to construct a taut
triangulation of $(M, \gamma)$.

Consider first the
case where $S$ has positive triangular number. Then, $S - \delta$ has an
ideal triangulation. We extend this to an ideal triangulation
$T$ of ${\cal R}_-(\gamma_S)$ and ${\cal R}_+(\gamma_S)$, which
is possible since each component of ${\cal R}_\pm(\gamma)
- {\rm int}({\cal N}(S))$ has positive triangular number. However, 
${\cal R}_-(\gamma_S)$ and ${\cal R}_+(\gamma_S)$ already
come equipped with an ideal triangulation, inherited from
the taut triangulation of $(M_S, \gamma_S)$.
Using elementary moves, we may alter this to $T$. Then
glue the two copies of $S$ in ${\cal R}_-(\gamma_S)$ and
${\cal R}_+(\gamma_S)$. We claim that the result is a taut triangulation of
$(M, \gamma)$. By construction, we have guaranteed
that each ideal 3--simplex has two inward-pointing faces
and two outward-pointing faces. Also, the transverse
orientations on the ideal 2--simplices in ${\cal R}_\pm(\gamma)$ are
as they should be. We now check that the orientations
of the faces around edges of $M$ are correct.

Consider an edge $e$ of $M$. 
If $e$ is in the interior of $M$, then there are
two possibilities: either it comes from a single edge
in the interior of $M_S$, in which case the faces around
it are already correctly oriented; or $e$ arises by identifying
two edges in $\partial M_S - \gamma_S$, one in ${\cal R}_-(\gamma_S)$
and one in ${\cal R}_+(\gamma_S)$, and so in this case
also, all but precisely two pairs of adjacent faces
around $e$ have cusped intersection. If $e$ is
an edge in $\partial M - \gamma$, then it came from
an edge in ${\cal R}_\pm(\gamma_S)$ and possibly
several edges of $\gamma_S - \delta_S$. Again, the
faces around $e$ are correctly oriented. Finally,
if $e$ is in $\gamma - \delta$, then it is formed
from one or more edges of $\gamma_S - \delta_S$,
and the faces around $e$ in this case all have cusped intersection. This
verifies that this is a taut
triangulation of $(M, \gamma)$.

The case where $S$ has zero triangular number is similar. 
Since $M$ has no $\delta$--essential annuli, $S$ is not an
annulus with its boundary in $\delta$. Hence, $S$
must be a disc intersecting $\delta$ twice. Since
$\gamma$ separates $\partial M - T(\gamma)$ into ${\cal R}_-(\gamma)$
and ${\cal R}_+(\gamma)$, either both arcs of
$\partial S \cap \delta$ are glancing or they
are both transverse. However, in the former case,
$S$ would be boundary-parallel, which is a contradiction.
Hence, both arcs are transverse and $S$ is a product disc. There are two copies of this
product disc in ${\cal R}_+(M_S)$ and ${\cal R}_-(M_S)$,
adjacent to two arcs of $\gamma_S - \delta_S$. Glue
these two edges together, forming an edge $e$. The result is not
quite a copy of $(M, \gamma)$, since it is not a 3--manifold in
a neighbourhood of $e$. Let $T_1$ and
$T_2$ be the two ideal triangles in ${\cal R}_-(M)$
adjacent to $e$. If necessary, we may perform a elementary
move to ${\cal R}_-(M)$ to ensure that $T_1$ and
$T_2$ are distinct. Then, attach an ideal 3--simplex
to $T_1$ and $T_2$, realizing an elementary move.
The result is now a taut triangulation
of $(M, \gamma)$. \endprf

We finish this section with a simple result.
When assigning a transverse
orientation to an ideal triangulation, one need not check
all the conditions of the definition of tautness; instead one can apply the
following proposition.

\proclaim{Proposition 9}Suppose that a transverse orientation
is assigned to the faces of an ideal triangulation of
a 3--manifold $M$, with $\partial M$ a collection of tori.
Then this specifies a taut ideal triangulation, providing that:
\items \itemb for each ideal 3--simplex, precisely two of its
faces are oriented into the 3--simplex, and precisely two
are oriented outwards, and
\itemb around each edge, at least one pair of
adjacent faces encircling the edge do not have cusped intersection.
\enditems\endproc

\prf The ideal triangulation of $M$ induces
a triangulation of $\partial M$. The transverse orientation on the
2--simplices of $M$ induces a transverse orientation on
the 1--simplices of $\partial M$. This specifies interior
angles of either $0$ or $\pi$ at each corner of each
triangle of $\partial M$. Let $V$, $E$ and $F$ be the
number of vertices, edges and faces of $\partial M$.
Let $N$ be the number of $\pi$ interior angles.
The first of the conditions in
the proposition gives that each triangle contains precisely
two zero interior angles and one $\pi$ interior angle, and
hence $N  = F$. 
Note that for orientation reasons, each vertex of $\partial M$
has an even number of $\pi$ interior angles.
Hence, the second of the above conditions gives that
each vertex of $\partial M$ has at least two
interior angles of $\pi$, and so $N \geq 2V$.
However, $0 = V - E + F = V - F/2$, and so $F = 2V$.
There must therefore be precisely two interior angles
of $\pi$ at each vertex of $\partial M$, and
hence this is a taut ideal triangulation. \endprf

\section{Singular surfaces and taut ideal triangulations}

In this section, we will prove that the singular genus of a
knot in $S^3$ is equal to its genus. We emphasise that
no further sutured manifold theory will be required. 
We now recall the outline of Gabai's original proof
of this result. More generally, he
showed that, for any compact orientable irreducible 3--manifold $M$
with $\partial M$ a (possibly empty) collection of tori,
the minimal Thurston complexity of embedded surfaces representing
a class in $H_2(M, \partial M)$ is equal to the minimal
Thurston complexity of singular surfaces representing that class.
He first constructed, for each non-zero class in $H_2(M,\partial M)$,
a taut transversely oriented foliation with a compact leaf $S$ (or leaves) representing
that class. This foliation defines a 2--plane field, which
has an associated Euler class in $e \in H^2(M, \partial M)$. 

Gabai then showed using a doubling argument
that it sufficed to consider the case
where $\partial M = \emptyset$. Then he considered how another
closed oriented surface $F$ mapped into $M$ 
interacts with the taut foliation,
where no component of $F$ is a 2--sphere and $F$ is {\sl homotopically incompressible}, 
meaning that the  only simple closed curves in $F$ which are 
homotopically trivial in $M$ are those which bound discs in $F$.
He showed that $F$ may be homotoped so that the non-transverse
intersections between $F$ and the foliation are of two
types: saddle tangencies and circle tangencies. The evaluation of
the cohomology class $e$ on the oriented surface $F$
is equal to the number of saddle singularities, counted
with sign, the sign depending on whether the transverse orientations of
$F$ and the foliation agree or disagree at a particular
saddle. If $[F,\partial F] = [S, \partial S] \in H_2(M, \partial M)$, then
of course $|e([F, \partial F])| = |e([S, \partial S])|$, which
is precisely $-\chi(S)$. Hence,
$-\chi(S)$ is the number of saddles of $F$ counted with sign.
However, $-\chi(F)$ is the number of saddles of $F$ 
counted without sign. Hence, $\chi_-(F) = -\chi(F) \geq -\chi(S)
= \chi_-(S)$.

Our combinatorial substitute for the taut foliation in 
the above argument is the taut ideal
triangulation of Theorem 2. We therefore need to consider how another
(possibly non-embedded) surface $F$ interacts with
this ideal triangulation. The obvious way to analyse this
is to use a version of normal surface theory. Similar considerations
arose when dealing with angled ideal triangulations in
[\ref{4}]. We recall the main points made there.

Since we are considering surfaces with boundary, it
is best to truncate each of the tetrahedra of the
ideal triangulation. The boundary of each truncated
tetrahedron $\Delta$ is decomposed into four triangles (which
are the intersection with $\partial M$) and four
hexagons (the truncated 2--simplices). Note that
$\partial \Delta$ inherits a 1--complex, which is
the union of the truncated 1--simplices and
the boundary of the triangles $\Delta \cap \partial M$.

Let $F$ be a compact orientable surface mapped into $M$, with 
$\partial F$ sent to $\partial M$.
Suppose that $F$ is homotopically incompressible and also
{\sl homotopically $\partial$--incompressible},
meaning that no embedded essential arc in $F$ can be homotoped 
in $M$ (keeping its endpoints fixed) to an arc in $\partial M$. 
Suppose also that $F$ contains no 2--sphere components
and no discs parallel to discs in $\partial M$.
We will see later (in Proposition 10) that $M$ is
irreducible and has incompressible boundary.
Hence, by the discussion in [\ref{4}], there is a homotopy of $F$ taking it
into {\sl admissible form}, which means that it
satisfies the following conditions for each truncated
3--simplex $\Delta$:
\items 
\item{(i)} $F \cap \Delta$ is a collection of 
discs in $\Delta$ intersecting $\partial \Delta$ in closed curves
transverse to the interior of the 1--cells in $\partial \Delta$;
\item{(ii)} no curve of $F \cap \partial \Delta$ is disjoint
from the 1--cells in $\partial \Delta$;
\item{(iii)} no arc of intersection between $F$ and 
a hexagon $H$ has endpoints lying in the same 1--cell
of $\partial \Delta$ or in adjacent 1--cells
of $\partial \Delta$; and
\item{(iv)} no arc of intersection between $F$ 
and a triangle of $\partial M \cap \Delta$
has endpoints lying in the same 1--cell
of $\partial \Delta$.
\enditems

\figure %%8
{\epsfbox{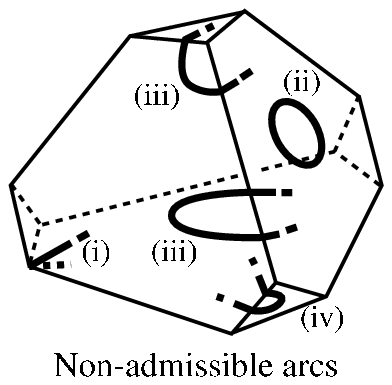}}
\endfigure

When $F$ is an admissible surface in a taut ideal triangulation, its intersection with
the 2--skeleton of $M - \partial M$ forms a transversely oriented
branched 1--manifold in $F$. Its complimentary
regions are discs. The boundary of such a disc $D$ 
inherits a number of {\sl cusps} which arise in
two possible situations:
\items \itemb either when $\partial D$ runs over a cusped
intersection between adjacent 2--simplices of $M$, or
\itemb when $\partial D$ runs over $\partial M$.
\enditems

If $c(D)$ is the number of cusps of $D$, we define 
the {\sl combinatorial area} of $D$ to be 
$${\rm Area}(D) = \pi(c(D) - 2).$$
This concurs with the definition of combinatorial area
in [\ref{4}], which was given in terms of interior angles.
(Recall that the transverse orientation on the ideal triangulation specifies an interior
angle of either zero or $\pi$ at each edge of each ideal
3--simplex.) It is clear that ${\rm Area}(D)$ is non-negative
for any admissible disc $D$.
We define ${\rm Area}(F)$ to be the sum of the
combinatorial areas of its discs. With this definition,
the argument in Proposition 4.3 of [\ref{4}] gives that
$${\rm Area}(F) = -2 \pi \chi(F).$$
Thus, we have the following immediate corollary.

\proclaim{Proposition 10}A taut ideal triangulation contains no admissible
2--spheres or discs. Hence, the underlying 3--manifold
is irreducible and its boundary is a collection of incompressible tori.
\endproc

Placing $F$ into admissible form is the analogue of Gabai's
method of homotoping $F$ so that its non-transverse intersections
with a foliation are saddles and centres. The above formula for
${\rm Area}(F)$ plays the r\^ole of the formula for the Euler characteristic
of $F$ in terms of the number of saddles.

We now need an analogue of the Euler class of a foliation.
Instead of finding a class in $H^2(M, \partial M)$, we
construct a class in $H_1(M)$. Let $G$ be the 4--valent
graph which is the 1--skeleton of the spine dual to the
ideal triangulation. The transverse orientation on
the ideal 2--simplices determines an orientation on
each edge of $G$. Since two edges point into
each vertex of $G$ and two edges point out,
this forms a 1--cycle $[G] \in H_1(M)$.
When $S$ is a surface as in Theorem 3 that is carried
by the underlying branched surface of the taut
ideal triangulation, then $G$ intersects each ideal 2--simplex of
$S - \partial S$ precisely once, and these
intersection points all have the same sign. Hence,
we have the following formula for the intersection
at the level of homology:
$$|[G] \cdot [S, \partial S]|= t(S, \partial S) = - 2 \chi(S).$$

The above formula and the following proposition will complete the proof
of Theorem 3.

\proclaim{Proposition 11}Let $F$ be a compact
orientable surface with no sphere or disc components.
Then for any map $(F, \partial F) \rightarrow (M, \partial M)$,
we have
$$|[G] \cdot [F, \partial F]| \leq - 2\chi(F).$$
\endproc

We may reduce to the case where $F$ is connected,
and where $[F, \partial F] \not= 0 \in H_2(M, \partial M)$.
If we homotopically compress and homotopically
$\partial$--compress $F$, this does not change its
class in $H_2(M, \partial M)$ and it increases its
Euler characteristic. Also, it does not create any
spheres or discs, since $M$ is irreducible and has
incompressible boundary, and $[F, \partial F] \not=0
\in H_2(M,\partial M)$. Hence, we may assume that
$F$ is homotopically incompressible and homotopically
$\partial$--incompressible. We may therefore homotope
it into admissible form. We may also homotope $F$
so that each point of intersection with each
1--cell of $\partial \Delta$ is at the midpoint
of the 1--cell, and each arc of intersection
between $F$ and a 2--cell of $\partial \Delta$ runs linearly between
these points. Hence, each truncated hexagonal 2--simplex
intersects $F$ only in one of nine
possible curves (shown in Figure 10), and $F$ intersects each triangle in $\partial M$
in at most three possible curves. Note that, for
each hexagon $H$, three of the possible curves of
$F \cap H$ run straight through $G \cap H$. This
causes a few minor technical problems. If $D$ is
a disc of $F \cap \Delta$ with $\partial D$ avoiding $G$,
then there is a well-defined signed intersection number
$G \cdot D$ between $G$ and $D$ which is
invariant under homotopies of $D$ in $\Delta$
that keep $\partial D$ fixed. If $\partial D$ hits $G$ in
a number of points, we may perturb $\partial D$ at each of these
points so as to miss $G$, in one of
two possible ways. After this perturbation, there
is then a well-defined intersection
number with $G$. We define $G \cdot D$ to be the
average  of these signed intersection numbers
over all such perturbations and all points of $D \cap G$.
It follows fairly rapidly from this definition that
$[G] \cdot [F,\partial F]$ is simply
the sum of $G \cdot D$ over all discs $D$ of $F \cap \Delta$
and all truncated 3--simplices $\Delta$ of $M$. Hence,
Proposition 11 follows from the following result.

\proclaim{Proposition 12}For any admissible disc
$D$ in $\Delta$, ${\rm Area}(D) \geq \pi |G \cdot D|$.\endproc

For this implies that
$$
|[G] \cdot [F, \partial F]| =
\left\vert \sum_D G \cdot D \right\vert
\leq \sum_D | G \cdot D |$$
$$\leq \sum_D {\rm Area}(D)/\pi
= {\rm Area}(F)/\pi = -2 \chi(F).$$

{\bf Proof of Proposition 12}\qua It is
possible to compute $G \cdot D$ in terms of
the arcs of intersection between $\partial D$
and the hexagons of $\partial \Delta$.
Label the four edges of $\Delta$
having zero interior angle with $e_1, \dots, e_4$.
For $i = 1$ to $4$, let $\alpha_i$ be the arc in $\partial \Delta$
running from a point of $G \cap \partial \Delta$ linearly to
the midpoint of $e_i$ and then continuing linearly on to
another point of $G \cap \partial \Delta$. Let $G_i$
be the arc in $G$ running between the endpoints of
$\alpha_i$. Orient $G_i$ according to the transverse orientation of the
2--simplices of $\Delta$. 

\figure %%9
{\epsfbox{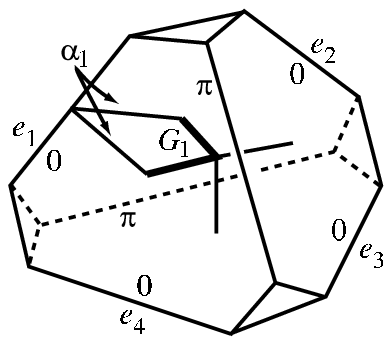}}
\endfigure

Note that
$$\sum_{i=1}^4 [G_i, \partial G_i] = 2 [G, G \cap \partial \Delta]
\in H_1(\Delta, G \cap \partial \Delta),$$
and so
$$\sum_{i=1}^4 G_i \cdot D = 2 (G \cdot D).$$
But $G_i \cdot D$ is (providing $\partial D$ misses $G_i \cap
\partial \Delta$) simply the winding number of
$\partial D$ around the annulus $\partial \Delta - G_i$,
which is the signed intersection number between
$\partial D$ and $\alpha_i$.
Hence, we can calculate $G \cdot D$ by computing, for
each arc $A$ of intersection between a hexagon $H$ of $\partial \Delta$
and $\partial D$, the signed intersection number between
$A$ and $\alpha_i \cap H$ (weighted by 1/2 if $A$ and $\alpha_i \cap H$
intersect at an endpoint of $\alpha_i \cap H$), and then summing these contributions
over all the $\alpha_i$, all the arcs $A$ and all the hexagons $H$,
and then dividing by two (since $\sum_i G_i \cdot D = 2 G \cdot D$).
Figure 10 shows the nine possibilities
for $A$ in each hexagon and the contribution that
each makes to $G \cdot D$:

\figure %%10
{\epsfbox{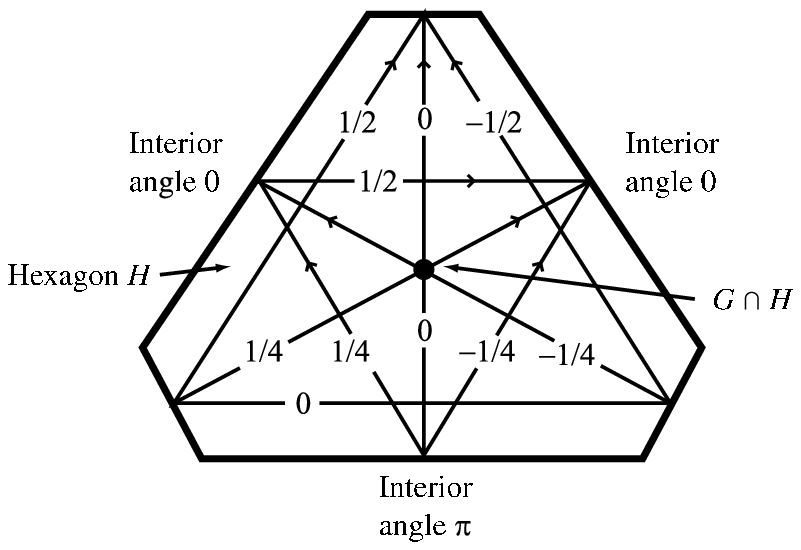}}
\endfigure

Note that the modulus of the contribution of each arc $A$ to $G \cdot D$ is
at most one quarter the number of cusps at the endpoints of $A$.
Hence, 
$$|G \cdot D| \leq c(D)/2.$$
When $c(D) \geq 4$, this proves the proposition, since then
$${\rm Area}(D) = \pi(c(D) - 2) \geq \pi c(D)/2 \geq \pi |G \cdot D|.$$
There are only a finitely many admissible discs $D$
with $c(D) < 4$. These are shown in Figure 11 (up to obvious
symmetries of $\Delta$) and are easily checked to
satisfy the proposition. \endprf

\vglue 15pt plus 20pt

\centerline{\epsfxsize 4.6truein\epsfbox{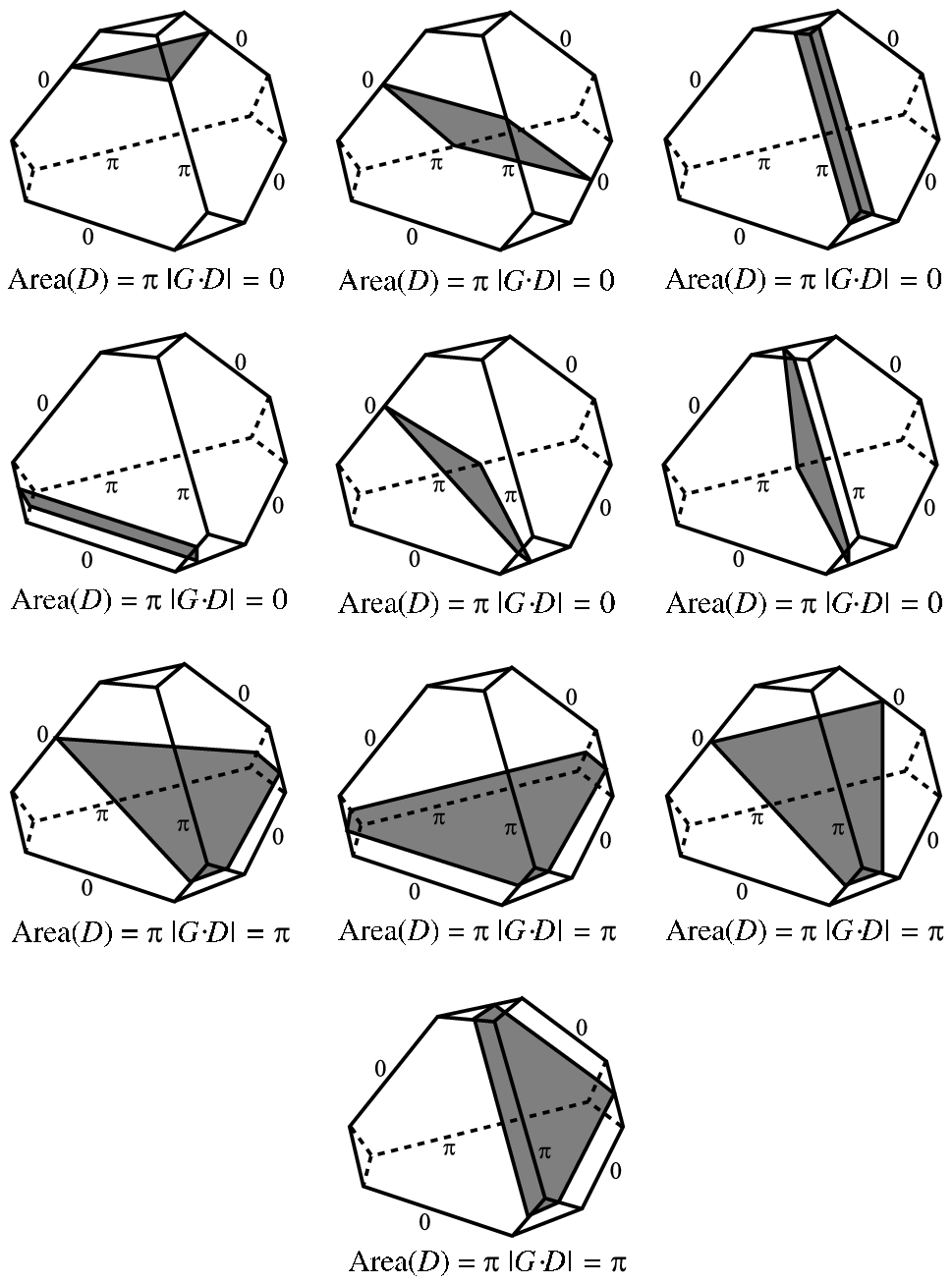}}
\centerline{\figurelabel}

\section{Further questions}

Although this paper presents simplified proofs of several
results that had previously been proved using foliation theory,
it is in no way meant as a replacement for that theory.
One of the principal limitations of taut ideal triangulations
is that they do not occur in closed 3--manifolds, whereas
taut foliations may of course arise. The first of the
following questions addresses this issue. The remaining
questions relate to other possible applications
of taut ideal triangulations.

(1)\qua Is there a version of taut ideal triangulations for
closed 3--manifolds? One candidate is the structures
on triangulations defined by Calegari in [\ref{1}].

(2)\qua Let $T$ be a taut ideal triangulation of a 3--manifold $M$,
with $\partial M$ a single torus. We say that a slope $s$ on
$\partial M$ is {\sl carried} by $T$ if there is a lamination
fully carried by the underlying branched surface of $T$
which intersects $\partial M$ only in simple closed curves
of slope $s$. Which slopes are carried by $T$? Is the set of slopes
open? Certainly, if a slope is carried by $T$, then
the manifold obtained by Dehn filling $M$ along that
slope has a taut transversely oriented foliation transverse to
the surgery curve. This yields topological restraints
on the possible slopes carried by any taut ideal
triangulation. For example, the meridian of a knot exterior
in $S^3$ is never carried by a taut ideal triangulation.

(3)\qua Let $M$ be a compact orientable irreducible atoroidal 
3--manifold with $\partial M$ a non-empty union of incompressible tori.
Does $M$ have a taut ideal triangulation whose
angles can be perturbed to give an angled ideal triangulation,
as in [\ref{4}]? Certainly, atoroidality is essential here.
A positive answer to this question would give a new
construction of angled ideal triangulations, which
might provide a useful insight into the geometrisation
of $M$.

\references

\bye